\catcode`@=11
\def\@height{height}
\def\@depth{depth}
\def\@width{width}

\newcount\@tempcnta
\newcount\@tempcntb

\newdimen\@tempdima
\newdimen\@tempdimb

\newbox\@tempboxa

\def\@ifnextchar#1#2#3{\let\@tempe #1\def\@tempa{#2}\def\@tempb{#3}\futurelet
    \@tempc\@ifnch}
\def\@ifnch{\ifx \@tempc \@sptoken \let\@tempd\@xifnch
      \else \ifx \@tempc \@tempe\let\@tempd\@tempa\else\let\@tempd\@tempb\fi
      \fi \@tempd}
\def\@ifstar#1#2{\@ifnextchar *{\def\@tempa*{#1}\@tempa}{#2}}

\def\@whilenoop#1{}
\def\@whilenum#1\do #2{\ifnum #1\relax #2\relax\@iwhilenum{#1\relax 
     #2\relax}\fi}
\def\@iwhilenum#1{\ifnum #1\let\@nextwhile=\@iwhilenum 
         \else\let\@nextwhile=\@whilenoop\fi\@nextwhile{#1}}

\def\@whiledim#1\do #2{\ifdim #1\relax#2\@iwhiledim{#1\relax#2}\fi}
\def\@iwhiledim#1{\ifdim #1\let\@nextwhile=\@iwhiledim 
        \else\let\@nextwhile=\@whilenoop\fi\@nextwhile{#1}}

\newdimen\@wholewidth
\newdimen\@halfwidth
\newdimen\unitlength \unitlength =1pt
\newbox\@picbox
\newdimen\@picht

\def\@nnil{\@nil}
\def\@empty{}
\def\@fornoop#1\@@#2#3{}

\def\@for#1:=#2\do#3{\edef\@fortmp{#2}\ifx\@fortmp\@empty \else
    \expandafter\@forloop#2,\@nil,\@nil\@@#1{#3}\fi}

\def\@forloop#1,#2,#3\@@#4#5{\def#4{#1}\ifx #4\@nnil \else
       #5\def#4{#2}\ifx #4\@nnil \else#5\@iforloop #3\@@#4{#5}\fi\fi}

\def\@iforloop#1,#2\@@#3#4{\def#3{#1}\ifx #3\@nnil 
       \let\@nextwhile=\@fornoop \else
      #4\relax\let\@nextwhile=\@iforloop\fi\@nextwhile#2\@@#3{#4}}

\def\@tfor#1:=#2\do#3{\xdef\@fortmp{#2}\ifx\@fortmp\@empty \else
    \@tforloop#2\@nil\@nil\@@#1{#3}\fi}
\def\@tforloop#1#2\@@#3#4{\def#3{#1}\ifx #3\@nnil 
       \let\@nextwhile=\@fornoop \else
      #4\relax\let\@nextwhile=\@tforloop\fi\@nextwhile#2\@@#3{#4}}

\def\@makepicbox(#1,#2){\leavevmode\@ifnextchar 
   [{\@imakepicbox(#1,#2)}{\@imakepicbox(#1,#2)[]}}

\long\def\@imakepicbox(#1,#2)[#3]#4{\vbox to#2\unitlength
   {\let\mb@b\vss \let\mb@l\hss\let\mb@r\hss
    \let\mb@t\vss
    \@tfor\@tempa :=#3\do{\expandafter\let
        \csname mb@\@tempa\endcsname\relax}%
\mb@t\hbox to #1\unitlength{\mb@l #4\mb@r}\mb@b}}

\def\picture(#1,#2){\@ifnextchar({\@picture(#1,#2)}{\@picture(#1,#2)(0,0)}}

\def\@picture(#1,#2)(#3,#4){\@picht #2\unitlength
\setbox\@picbox\hbox to #1\unitlength\bgroup 
\hskip -#3\unitlength \lower #4\unitlength \hbox\bgroup\ignorespaces}

\def\endpicture{\egroup\hss\egroup\ht\@picbox\@picht
\dp\@picbox\z@\leavevmode\box\@picbox}

\long\def\put(#1,#2)#3{\@killglue\raise#2\unitlength\hbox to \z@{\kern
#1\unitlength #3\hss}\ignorespaces}

\long\def\multiput(#1,#2)(#3,#4)#5#6{\@killglue\@multicnt=#5\relax
\@xdim=#1\unitlength
\@ydim=#2\unitlength
\@whilenum \@multicnt > 0\do
{\raise\@ydim\hbox to \z@{\kern
\@xdim #6\hss}\advance\@multicnt \m@ne\advance\@xdim
#3\unitlength\advance\@ydim #4\unitlength}\ignorespaces}

\def\@killglue{\unskip\@whiledim \lastskip >\z@\do{\unskip}}

\def\thinlines{\let\@linefnt\tenln \let\@circlefnt\tencirc
  \@wholewidth\fontdimen8\tenln \@halfwidth .5\@wholewidth}
\def\thicklines{\let\@linefnt\tenlnw \let\@circlefnt\tencircw
  \@wholewidth\fontdimen8\tenlnw \@halfwidth .5\@wholewidth}

\def\linethickness#1{\@wholewidth #1\relax \@halfwidth .5\@wholewidth}

\def\shortstack{\@ifnextchar[{\@shortstack}{\@shortstack[c]}}

\def\@shortstack[#1]{\leavevmode
\vbox\bgroup\baselineskip-1pt\lineskip 3pt\let\mb@l\hss
\let\mb@r\hss \expandafter\let\csname mb@#1\endcsname\relax
\let\\\@stackcr\@ishortstack}

\def\@ishortstack#1{\halign{\mb@l ##\unskip\mb@r\cr #1\crcr}\egroup}

\def\@stackcr{\@ifstar{\@ixstackcr}{\@ixstackcr}}
\def\@ixstackcr{\@ifnextchar[{\@istackcr}{\cr\ignorespaces}}

\def\@istackcr[#1]{\cr\noalign{\vskip #1}\ignorespaces}

\newif\if@negarg

\def\droite(#1,#2)#3{\@xarg #1\relax \@yarg #2\relax
\@linelen=#3\unitlength
\ifnum\@xarg =0 \@vline 
  \else \ifnum\@yarg =0 \@hline \else \@sline\fi
\fi}

\def\@sline{\ifnum\@xarg< 0 \@negargtrue \@xarg -\@xarg \@yyarg -\@yarg
  \else \@negargfalse \@yyarg \@yarg \fi
\ifnum \@yyarg >0 \@tempcnta\@yyarg \else \@tempcnta -\@yyarg \fi
\ifnum\@tempcnta>6 \@badlinearg\@tempcnta0 \fi
\ifnum\@xarg>6 \@badlinearg\@xarg 1 \fi
\setbox\@linechar\hbox{\@linefnt\@getlinechar(\@xarg,\@yyarg)}%
\ifnum \@yarg >0 \let\@upordown\raise \@clnht\z@
   \else\let\@upordown\lower \@clnht \ht\@linechar\fi
\@clnwd=\wd\@linechar
\if@negarg \hskip -\wd\@linechar \def\@tempa{\hskip -2\wd\@linechar}\else
     \let\@tempa\relax \fi
\@whiledim \@clnwd <\@linelen \do
  {\@upordown\@clnht\copy\@linechar
   \@tempa
   \advance\@clnht \ht\@linechar
   \advance\@clnwd \wd\@linechar}%
\advance\@clnht -\ht\@linechar
\advance\@clnwd -\wd\@linechar
\@tempdima\@linelen\advance\@tempdima -\@clnwd
\@tempdimb\@tempdima\advance\@tempdimb -\wd\@linechar
\if@negarg \hskip -\@tempdimb \else \hskip \@tempdimb \fi
\multiply\@tempdima \@m
\@tempcnta \@tempdima \@tempdima \wd\@linechar \divide\@tempcnta \@tempdima
\@tempdima \ht\@linechar \multiply\@tempdima \@tempcnta
\divide\@tempdima \@m
\advance\@clnht \@tempdima
\ifdim \@linelen <\wd\@linechar
   \hskip \wd\@linechar
  \else\@upordown\@clnht\copy\@linechar\fi}

\def\@hline{\ifnum \@xarg <0 \hskip -\@linelen \fi
\vrule \@height \@halfwidth \@depth \@halfwidth \@width \@linelen
\ifnum \@xarg <0 \hskip -\@linelen \fi}

\def\@getlinechar(#1,#2){\@tempcnta#1\relax\multiply\@tempcnta 8
\advance\@tempcnta -9 \ifnum #2>0 \advance\@tempcnta #2\relax\else
\advance\@tempcnta -#2\relax\advance\@tempcnta 64 \fi
\char\@tempcnta}

\def\vector(#1,#2)#3{\@xarg #1\relax \@yarg #2\relax
\@tempcnta \ifnum\@xarg<0 -\@xarg\else\@xarg\fi
\ifnum\@tempcnta<5\relax
\@linelen=#3\unitlength
\ifnum\@xarg =0 \@vvector 
  \else \ifnum\@yarg =0 \@hvector \else \@svector\fi
\fi
\else\@badlinearg\fi}

\def\@hvector{\@hline\hbox to 0pt{\@linefnt 
\ifnum \@xarg <0 \@getlarrow(1,0)\hss\else
    \hss\@getrarrow(1,0)\fi}}

\def\@vvector{\ifnum \@yarg <0 \@downvector \else \@upvector \fi}

\def\@svector{\@sline
\@tempcnta\@yarg \ifnum\@tempcnta <0 \@tempcnta=-\@tempcnta\fi
\ifnum\@tempcnta <5
  \hskip -\wd\@linechar
  \@upordown\@clnht \hbox{\@linefnt  \if@negarg 
  \@getlarrow(\@xarg,\@yyarg) \else \@getrarrow(\@xarg,\@yyarg) \fi}%
\else\@badlinearg\fi}

\def\@getlarrow(#1,#2){\ifnum #2 =\z@ \@tempcnta='33\else
\@tempcnta=#1\relax\multiply\@tempcnta \sixt@@n \advance\@tempcnta
-9 \@tempcntb=#2\relax\multiply\@tempcntb \tw@
\ifnum \@tempcntb >0 \advance\@tempcnta \@tempcntb\relax
\else\advance\@tempcnta -\@tempcntb\advance\@tempcnta 64
\fi\fi\char\@tempcnta}

\def\@getrarrow(#1,#2){\@tempcntb=#2\relax
\ifnum\@tempcntb < 0 \@tempcntb=-\@tempcntb\relax\fi
\ifcase \@tempcntb\relax \@tempcnta='55 \or 
\ifnum #1<3 \@tempcnta=#1\relax\multiply\@tempcnta
24 \advance\@tempcnta -6 \else \ifnum #1=3 \@tempcnta=49
\else\@tempcnta=58 \fi\fi\or 
\ifnum #1<3 \@tempcnta=#1\relax\multiply\@tempcnta
24 \advance\@tempcnta -3 \else \@tempcnta=51\fi\or 
\@tempcnta=#1\relax\multiply\@tempcnta
\sixt@@n \advance\@tempcnta -\tw@ \else
\@tempcnta=#1\relax\multiply\@tempcnta
\sixt@@n \advance\@tempcnta 7 \fi\ifnum #2<0 \advance\@tempcnta 64 \fi
\char\@tempcnta}

\def\@vline{\ifnum \@yarg <0 \@downline \else \@upline\fi}

\def\@upline{\hbox to \z@{\hskip -\@halfwidth \vrule \@width \@wholewidth
   \@height \@linelen \@depth \z@\hss}}

\def\@downline{\hbox to \z@{\hskip -\@halfwidth \vrule \@width \@wholewidth
   \@height \z@ \@depth \@linelen \hss}}

\def\@upvector{\@upline\setbox\@tempboxa\hbox{\@linefnt\char'66}\raise 
     \@linelen \hbox to\z@{\lower \ht\@tempboxa\box\@tempboxa\hss}}

\def\@downvector{\@downline\lower \@linelen
      \hbox to \z@{\@linefnt\char'77\hss}}

\def\dashbox#1(#2,#3){\leavevmode\hbox to \z@{\baselineskip \z@%
\lineskip \z@%
\@dashdim=#2\unitlength%
\@dashcnt=\@dashdim \advance\@dashcnt 200
\@dashdim=#1\unitlength\divide\@dashcnt \@dashdim
\ifodd\@dashcnt\@dashdim=\z@%
\advance\@dashcnt \@ne \divide\@dashcnt \tw@ 
\else \divide\@dashdim \tw@ \divide\@dashcnt \tw@
\advance\@dashcnt \m@ne
\setbox\@dashbox=\hbox{\vrule \@height \@halfwidth \@depth \@halfwidth
\@width \@dashdim}\put(0,0){\copy\@dashbox}%
\put(0,#3){\copy\@dashbox}%
\put(#2,0){\hskip-\@dashdim\copy\@dashbox}%
\put(#2,#3){\hskip-\@dashdim\box\@dashbox}%
\multiply\@dashdim 3 
\fi
\setbox\@dashbox=\hbox{\vrule \@height \@halfwidth \@depth \@halfwidth
\@width #1\unitlength\hskip #1\unitlength}\@tempcnta=0
\put(0,0){\hskip\@dashdim \@whilenum \@tempcnta <\@dashcnt
\do{\copy\@dashbox\advance\@tempcnta \@ne }}\@tempcnta=0
\put(0,#3){\hskip\@dashdim \@whilenum \@tempcnta <\@dashcnt
\do{\copy\@dashbox\advance\@tempcnta \@ne }}%
\@dashdim=#3\unitlength%
\@dashcnt=\@dashdim \advance\@dashcnt 200
\@dashdim=#1\unitlength\divide\@dashcnt \@dashdim
\ifodd\@dashcnt \@dashdim=\z@%
\advance\@dashcnt \@ne \divide\@dashcnt \tw@
\else
\divide\@dashdim \tw@ \divide\@dashcnt \tw@
\advance\@dashcnt \m@ne
\setbox\@dashbox\hbox{\hskip -\@halfwidth
\vrule \@width \@wholewidth 
\@height \@dashdim}\put(0,0){\copy\@dashbox}%
\put(#2,0){\copy\@dashbox}%
\put(0,#3){\lower\@dashdim\copy\@dashbox}%
\put(#2,#3){\lower\@dashdim\copy\@dashbox}%
\multiply\@dashdim 3
\fi
\setbox\@dashbox\hbox{\vrule \@width \@wholewidth 
\@height #1\unitlength}\@tempcnta0
\put(0,0){\hskip -\@halfwidth \vbox{\@whilenum \@tempcnta < \@dashcnt
\do{\vskip #1\unitlength\copy\@dashbox\advance\@tempcnta \@ne }%
\vskip\@dashdim}}\@tempcnta0
\put(#2,0){\hskip -\@halfwidth \vbox{\@whilenum \@tempcnta< \@dashcnt
\relax\do{\vskip #1\unitlength\copy\@dashbox\advance\@tempcnta \@ne }%
\vskip\@dashdim}}}\@makepicbox(#2,#3)}

\newif\if@ovt 
\newif\if@ovb 
\newif\if@ovl 
\newif\if@ovr 
\newdimen\@ovxx
\newdimen\@ovyy
\newdimen\@ovdx
\newdimen\@ovdy
\newdimen\@ovro
\newdimen\@ovri

\def\@getcirc#1{\@tempdima #1\relax \advance\@tempdima 2pt\relax
  \@tempcnta\@tempdima
  \@tempdima 4pt\relax \divide\@tempcnta\@tempdima
  \ifnum \@tempcnta > 10\relax \@tempcnta 10\relax\fi
  \ifnum \@tempcnta >\z@ \advance\@tempcnta\m@ne
    \else \@warning{Oval too small}\fi
  \multiply\@tempcnta 4\relax
  \setbox \@tempboxa \hbox{\@circlefnt
  \char \@tempcnta}\@tempdima \wd \@tempboxa}

\def\@put#1#2#3{\raise #2\hbox to \z@{\hskip #1#3\hss}}

\def\oval(#1,#2){\@ifnextchar[{\@oval(#1,#2)}{\@oval(#1,#2)[]}}

\def\@oval(#1,#2)[#3]{\begingroup\boxmaxdepth \maxdimen
  \@ovttrue \@ovbtrue \@ovltrue \@ovrtrue
  \@tfor\@tempa :=#3\do{\csname @ov\@tempa false\endcsname}\@ovxx
  #1\unitlength \@ovyy #2\unitlength
  \@tempdimb \ifdim \@ovyy >\@ovxx \@ovxx\else \@ovyy \fi
  \advance \@tempdimb -2pt\relax  
  \@getcirc \@tempdimb
  \@ovro \ht\@tempboxa \@ovri \dp\@tempboxa
  \@ovdx\@ovxx \advance\@ovdx -\@tempdima \divide\@ovdx \tw@
  \@ovdy\@ovyy \advance\@ovdy -\@tempdima \divide\@ovdy \tw@
  \@circlefnt \setbox\@tempboxa
  \hbox{\if@ovr \@ovvert32\kern -\@tempdima \fi
  \if@ovl \kern \@ovxx \@ovvert01\kern -\@tempdima \kern -\@ovxx \fi
  \if@ovt \@ovhorz \kern -\@ovxx \fi
  \if@ovb \raise \@ovyy \@ovhorz \fi}\advance\@ovdx\@ovro
  \advance\@ovdy\@ovro \ht\@tempboxa\z@ \dp\@tempboxa\z@
  \@put{-\@ovdx}{-\@ovdy}{\box\@tempboxa}%
  \endgroup}

\def\@ovvert#1#2{\vbox to \@ovyy{%
    \if@ovb \@tempcntb \@tempcnta \advance \@tempcntb by #1\relax
      \kern -\@ovro \hbox{\char \@tempcntb}\nointerlineskip
    \else \kern \@ovri \kern \@ovdy \fi
    \leaders\vrule width \@wholewidth\vfil \nointerlineskip
    \if@ovt \@tempcntb \@tempcnta \advance \@tempcntb by #2\relax
      \hbox{\char \@tempcntb}%
    \else \kern \@ovdy \kern \@ovro \fi}}

\def\@ovhorz{\hbox to \@ovxx{\kern \@ovro
    \if@ovr \else \kern \@ovdx \fi
    \leaders \hrule height \@wholewidth \hfil
    \if@ovl \else \kern \@ovdx \fi
    \kern \@ovri}}

\def\circle{\@ifstar{\@dot}{\@circle}}
\def\@circle#1{\begingroup \boxmaxdepth \maxdimen \@tempdimb #1\unitlength
   \ifdim \@tempdimb >15.5pt\relax \@getcirc\@tempdimb
      \@ovro\ht\@tempboxa 
     \setbox\@tempboxa\hbox{\@circlefnt
      \advance\@tempcnta\tw@ \char \@tempcnta
      \advance\@tempcnta\m@ne \char \@tempcnta \kern -2\@tempdima
      \advance\@tempcnta\tw@
      \raise \@tempdima \hbox{\char\@tempcnta}\raise \@tempdima
        \box\@tempboxa}\ht\@tempboxa\z@ \dp\@tempboxa\z@
      \@put{-\@ovro}{-\@ovro}{\box\@tempboxa}%
   \else  \@circ\@tempdimb{96}\fi\endgroup}

\def\@dot#1{\@tempdimb #1\unitlength \@circ\@tempdimb{112}}

\def\@circ#1#2{\@tempdima #1\relax \advance\@tempdima .5pt\relax
   \@tempcnta\@tempdima \@tempdima 1pt\relax
   \divide\@tempcnta\@tempdima 
   \ifnum\@tempcnta > 15\relax \@tempcnta 15\relax \fi    
   \ifnum \@tempcnta >\z@ \advance\@tempcnta\m@ne\fi
   \advance\@tempcnta #2\relax
   \@circlefnt \char\@tempcnta}

\font\tenln line10
\font\tencirc lcircle10
\font\tenlnw linew10
\font\tencircw lcirclew10

\thinlines   

\newcount\@xarg
\newcount\@yarg
\newcount\@yyarg
\newcount\@multicnt 
\newdimen\@xdim
\newdimen\@ydim
\newbox\@linechar
\newdimen\@linelen
\newdimen\@clnwd
\newdimen\@clnht
\newdimen\@dashdim
\newbox\@dashbox
\newcount\@dashcnt
\catcode`@=12

\overfullrule=0pt
\magnification=1200
\hsize=11.25cm    
\vsize=18cm
\hoffset=1cm
\font\grand=cmr10 at 14pt
\font\petit=cmr10 at 8pt
\font\normal=cmr10 at 10pt

\def\N{\noindent}
\def\S{\smallskip \par}
\def\M{\medskip \par}
\def\B{\bigskip \par}
\def\BB{\bigskip\bigskip\par}

\def\l{\dashv }
\def\r{\vdash }
\def\u{\underline }
\def\b{\backslash }
\def\m{\, \bot\,  }
\def\YY{Y\!\!\!\! Y}
\def\TT{T\!\!\! T}
\def\SS{\Sigma}
\def\oo{\omega}

\def\NN{{\bf N}}

\def\sqr#1#2{{\vcenter{\vbox{\hrule height.#2pt
\hbox{\vrule width .#2pt height#1pt \kern#1pt
\vrule width.#2pt}
\hrule height.#2pt}}}}
\def\square{\mathchoice\sqr64\sqr64\sqr{4.2}3\sqr33}

\catcode`@=11
\font\@linefnt linew10 at 2.4pt
\catcode`@=12

\def\arbreun{\kern-0.4ex
\hbox{\unitlength=.25pt
\picture(60,40)(0,0)
\put(30,0){\droite(0,1){10}}
\put(30,10){\droite(-1,1){30}}
\put(30,10){\droite(1,1){30}}
\put(20,20){\droite(1,1){20}}
\put(10,30){\droite(1,1){10}}
\endpicture}\kern 0.4ex}

\def\arbredeux{\kern-0.4ex
\hbox{\unitlength=.25pt
\picture(60,40)(0,0)
\put(30,0){\droite(0,1){10}}
\put(30,10){\droite(-1,1){30}}
\put(30,10){\droite(1,1){30}}
\put(20,20){\droite(1,1){20}}
\put(30,30){\droite(-1,1){10}}
\endpicture}\kern 0.4ex}

\def\arbretrois{\kern-0.4ex
\hbox{\unitlength=.25pt
\picture(60,40)(0,0)
\put(30,0){\droite(0,1){10}}
\put(30,10){\droite(-1,1){30}}
\put(30,10){\droite(1,1){30}}
\put(50,30){\droite(-1,1){10}}
\put(10,30){\droite(1,1){10}}
\endpicture}\kern 0.4ex}

\def\arbrequatre{\kern-0.4ex
\hbox{\unitlength=.25pt
\picture(60,40)(0,0)
\put(30,0){\droite(0,1){10}}
\put(30,10){\droite(-1,1){30}}
\put(30,10){\droite(1,1){30}}
\put(40,20){\droite(-1,1){20}}
\put(30,30){\droite(1,1){10}}
\endpicture}\kern 0.4ex}

\def\arbrecinq{\kern-0.4ex
\hbox{\unitlength=.25pt
\picture(60,40)(0,0)
\put(30,0){\droite(0,1){10}}
\put(30,10){\droite(-1,1){30}}
\put(30,10){\droite(1,1){30}}
\put(40,20){\droite(-1,1){20}}
\put(50,30){\droite(-1,1){10}}
\endpicture}\kern 0.4ex}

\def\arbreA{\kern-0.4ex
\hbox{\unitlength=.25pt
\picture(60,40)(0,0)
\put(30,0){\droite(0,1){10}}
\put(30,10){\droite(-1,1){30}}
\put(30,10){\droite(1,1){30}}
\endpicture}\kern 0.4ex}

\def\arbreB{\kern-0.4ex
\hbox{\unitlength=.25pt
\picture(60,40)(0,0)
\put(30,0){\droite(0,1){10}}
\put(30,10){\droite(-1,1){30}}
\put(30,10){\droite(1,1){30}}
\put(15,25){\droite(1,1){15}}
\endpicture}\kern 0.4ex}

\def\arbreC{\kern-0.4ex
\hbox{\unitlength=.25pt
\picture(60,40)(0,0)
\put(30,0){\droite(0,1){10}}
\put(30,10){\droite(-1,1){30}}
\put(30,10){\droite(1,1){30}}
\put(45,25){\droite(-1,1){15}}
\endpicture}\kern 0.4ex}

\def\trre{\kern-0.4ex
\hbox{\unitlength=1pt
\picture(100,70)(0,0)
\put(50,0){\droite(0,1){20}}
\put(50,20){\droite(-1,1){50}}
\put(50,20){\droite(1,1){50}}
\put(30,40){\droite(1,1){30}}
\put(10,60){\droite(1,1){10}}
\put(50,60){\droite(-1,1){10}}
\put(90,60){\droite(-1,1){10}}
\endpicture}\kern 0.4ex}

\def\trree{\kern-0.4ex
\hbox{\unitlength=1pt
\picture(100,70)(0,0)
\put(50,0){\droite(0,1){20}}
\put(50,20){\droite(-1,1){50}}
\put(50,20){\droite(1,1){50}}
\put(30,40){\droite(1,1){30}}
\put(10,60){\droite(1,1){10}}
\put(30,40){\droite(0,1){30}}
\put(90,60){\droite(-1,1){10}}
\put(50,20){\droite(1,2){25}}

\endpicture}\kern 0.4ex}

\def\arbreAa{\kern-0.4ex
\hbox{\unitlength=.25pt
\picture(60,80)(0,0)
\put(-10,50){$x$}
\put(50,50){$y$}
\put(30,0){\droite(0,1){10}}
\put(30,10){\droite(-1,1){30}}
\put(30,10){\droite(1,1){30}}
\endpicture}\kern 0.4ex}

\def\arbreAaa{\kern-0.4ex
\hbox{\unitlength=.25pt
\picture(60,80)(0,0)
\put(-20,60){$x^{(0)}$}
\put(110,60){$x^{(k)}$}
\put(60,-10){\droite(0,1){20}}
\put(60,10){\droite(-2,1){50}}
\put(40,30){$\cdots$}
\put(60,10){\droite(2,1){50}}
\endpicture}\kern 0.4ex}

\def\arbreunsix{\kern-0.4ex
\hbox{\unitlength=.25pt
\picture(200,60)(0,0)
\put(80,0){\droite(0,1){20}}
\put(80,20){\droite(-1,1){80}}
\put(80,20){\droite(1,1){80}}
\put(50,50){\droite(1,1){50}}
\put(80,80){\droite(-1,1){20}}
\put(140,80){\droite(-1,1){20}}
\put(40,60){$x$}
\put(70,60){\droite(-1,1){10}}
\put(70,100){$y^l$}
\put(130,100){$y^r$}

\put(160,20){$\longrightarrow$}

\put(300,0){\droite(0,1){20}}
\put(300,20){\droite(-1,1){80}}
\put(300,20){\droite(1,1){80}}
\put(300,80){\droite(1,1){20}}
\put(330,50){\droite(-1,1){50}}
\put(360,80){\droite(-1,1){20}}
\put(290,30){$x$}
\put(320,30){\droite(-1,1){10}}
\put(290,100){$y^l$}
\put(350,100){$y^r$}
\endpicture}\kern 0.4ex}

\def\arbredeuxdeux{\kern-0.4ex
\hbox{\unitlength=.25pt
\picture(60,40)(0,0)
\put(30,0){\droite(0,1){40}}
\put(30,10){\droite(-1,1){30}}
\put(30,10){\droite(1,1){30}}
\endpicture}\kern 0.4ex}

\def\arbreuntroistrois{\kern-0.4ex
\hbox{\unitlength=.25pt
\picture(60,40)(0,0)
\put(30,0){\droite(0,1){40}}
\put(30,10){\droite(-1,1){30}}
\put(30,10){\droite(1,1){30}}
\put(10,30){\droite(1,1){10}}
\endpicture}\kern 0.4ex}

\def\arbretroistroistrois{\kern-0.4ex
\hbox{\unitlength=.25pt
\picture(60,40)(0,0)
\put(30,0){\droite(0,1){10}}
\put(30,10){\droite(-2,1){60}}
\put(30,10){\droite(-2,3){20}}
\put(30,10){\droite(2,3){20}}
\put(30,10){\droite(2,1){60}}
\endpicture}\kern 0.4ex}

\def\arbreneufsixun{\kern-0.4ex
\hbox{\unitlength=.25pt
\picture(1000,200)(0,0)
\put(-40,0){$t=$}
\put(120,0){\droite(0,1){40}}
\put(120,40){\droite(-1,1){120}}
\put(120,40){\droite(1,1){120}}
\put(180,100){\droite(-1,1){30}}
\put(100, 140){$\cdots$}

\put(300,0){$t+1=$}
\put(520,0){\droite(0,1){40}}
\put(520,40){\droite(-1,1){120}}
\put(520,40){\droite(1,1){120}}
\put(580,100){\droite(-1,1){30}}
\put(500,140){$\cdots$}
\put(610,130){\droite(2,1){60}}

\put(700,0){$\cup$}
\put(920,0){\droite(0,1){40}}
\put(920,40){\droite(-1,1){120}}
\put(920,40){\droite(1,1){120}}
\put(980,100){\droite(-1,1){30}}
\put(900,140){$\cdots$}
\put(980,100){\droite(2,1){120}}

\endpicture}\kern 0.4ex}

\def\arbreneufsixdeux{\kern-0.4ex
\hbox{\unitlength=.25pt
\picture(1000,200)(0,0)

\put(50,0){$\cup$}
\put(220,0){\droite(0,1){40}}
\put(220,40){\droite(-1,1){120}}
\put(220,40){\droite(1,1){120}}
\put(280,100){\droite(-1,1){30}}
\put(200, 140){$\cdots$}
\put(250,70){\droite(2,1){180}}

\put(450,0){$\cup$}
\put(620,0){\droite(0,1){40}}
\put(620,40){\droite(-1,1){120}}
\put(620,40){\droite(1,1){120}}
\put(680,100){\droite(-1,1){30}}
\put(600,140){$\cdots$}
\put(620,40){\droite(2,1){180}}

\put(850,0){$\cup$}
\put(1020,0){\droite(0,1){40}}
\put(1020,40){\droite(-1,1){120}}
\put(1020,40){\droite(1,1){120}}
\put(1080,100){\droite(-1,1){30}}
\put(1000,140){$\cdots$}
\put(1020,20){\droite(2,1){180}}

\endpicture}\kern 0.4ex}

{\it \hfill \'A Claudio Procesi}
\BB
\centerline {\grand Arithmetree}
\BB

\centerline {\bf Jean-Louis Loday}
\BB

\N {\bf Abstract.} {\it    We construct an addition and a multiplication on the set
 of planar binary trees, closely related  to addition and multiplication on the integers.
This gives rise to a new kind of (noncommutative) arithmetic theory. The price to pay for this generalization is
that, first the addition is not commutative, second the multiplication is distributive with the addition only on the
left. This algebraic structure is the ``exponent part" of the free dendriform algebra on one generator, a notion related to
several other types of algebras.  

In the second part we extend this theory to all the planar trees. Then it is  related to the free dendriform trialgebra as
constructed in [LR3].}
\M

\N {\bf Introduction.} Elementary arithmetic deals with the natural numbers: 
$$0, 1, 2, \cdots , n , \cdots, $$
 on which one knows how to define an addition + and a
mutiplication $\times$. In this paper we propose the following generalization: we replace the integers by the 
planar binary trees. Recall that there are  $c_n={(2n!)\over n!(n+1)!}$ planar binary trees with $n+1$ leaves. The
integer $c_n$ is classically called the Catalan number. We first construct the sum of two planar binary trees. In
general this sum is not just a tree but a union of planar binary trees. However it comes with the following feature:
all the trees appearing in this sum are different. In other words it is a subset of the set of planar binary trees. We call
such a subset  a {\it grove}. We show that this sum indeed extends to groves and is associative. It is not
commutative, but there is an involution compatible with the sum. The construction of the sum takes advantage of a poset
structure on the set of planar binary trees.

Next we show that there also exists a multiplication on planar binary trees. Again, the product of two trees is not a tree in
general, but a grove. We show that the product of two groves is still a grove. This product is associative, distributive on
the left with the sum but not distributive on the right. 

The existence of this multiplication is due to a very peculiar
property of the addition. The sum of two planar binary trees turns out to be in fact the union of the results of two other
operations. Roughly speaking it is like making a difference between adding
$x$ on the left to
$y$  and adding $y$ on the right to $x$ . These two operations happen to satisfy some relations. When we take the
polynomial algebra with planar binary trees as exponents in place of integers, then we get what we
call  a {\it dendriform algebra}. The fact that this dendriform algebra is nothing but the free dendriform algebra on one
generator enables us to define the multiplication on planar binary trees.

The set of integers is very often used as an indexing set. However sometimes it is not sufficient and one has to move to
planar binary trees. This happens for instance in solving differential equations by means of series (cf. [Br], [BF]), and in
algebraic topology (generalization of the simplicial category (cf. [Fr]), of operads, of PROPs). As soon as one wants to
manipulate these objects, one needs to add and multiply the planar binary trees. This is one of the motivations of the
present work.

It is tempting to find out whether the addition (resp. the multiplication) can easily be described on some
of the other interpretations of the Catalan sets. We give one of them by introducing the ``permutation-like notation" of the
elements in the Catalan sets.
\B
In the second part of this paper we extend this arithmetic to {\it all} the planar trees, the case of planar binary trees
becoming a quotient of it. Since there are 3 different planar trees with three leaves, this theory is related to a type of
algebra defined by 3 operations. They are called {\it dendriform trialgebras} and where introduced in [LR3]. 
\B

 Thanks to Patrick Ion for suggesting the
terminology ``grove". Though it is not apparent in the text, this paper owes much to the book ``On Numbers
and Games" by J.H.C. Conway [Co].
\B
\centerline {\bf I. Arithmetic of planar binary trees}
\B

\N {\bf 1. The poset of planar binary trees.} 
\M
\N {\bf 1.1. Catalan sets.} Let $Y_0$ be a set with one element. The sets $Y_n$ for
$n\geq 1$ are defined inductively by the formula
$$Y_n:=Y_{n-1}\times Y_0 \cup \cdots \cup Y_{n-i}\times Y_{i-1}  \cup \cdots \cup Y_0\times Y_{n-1}.$$
If we denote by $a$ the unique element of $Y_0$, then an element of $Y_n$ can be described as a
(complete) parenthesizing of the word $aa\cdots  a$ of length $n+1$. Let $x\in Y_p$ and $y\in Y_q$. The element
$(x,y)\in Y_p\times Y_q$ viewed as an element in $Y_{p+q+1}$ is denoted
$x\vee y \in Y_{p+q+1}$. In terms of parenthesizing it simply consists in concatenating the two words and putting
a parenthesis at both ends.

There are many other combinatorial descriptions of the sets $Y_n$. We will
use two of them as described below, one classical: the planar binary trees, and one less classical: the
permutation-like notation. Others include: the triangulations of an $(n+2)$-gon, the vertices of the Stasheff
polytope of dimension
$n-1$, see [St] for many more.

Let $c_n$ be the number of elements of $Y_n$. It comes immediately: $c_0$ =1 and 
$$c_n:=c_{n-1} c_0 + \cdots + c_{n-i}  c_{i-1}  + \cdots + c_0 c_{n-1}.\eqno (1.1.1)$$
Hence the generating series $f(x):= \sum_{n\geq 0} c_n x^n$ satisfies the functional equation $ xf(x)^2 = f(x) -1$,
and we get $f(x) = {1-\sqrt {1-4x}\over 2x}$. As a consequence we get $c_n= {(2n)!\over
n!(n+1)!}$. It is  classically called the Catalan number, so $Y_n$ is called the {\it Catalan set}. 
\M

\N {\bf 1.2. Planar binary trees.} A {\it planar binary tree} (p.b.tree for short) is an oriented planar
graph drawn in the plane with $n+1$ leaves and one root, such that each internal
vertex has two leaves and one root. We consider these trees up to planar isotopy.
Here is an example:
\M

$$\trre$$

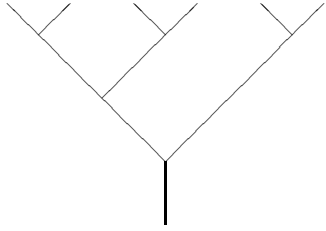
\centerline {figure 1}
\M
The number $n$ of internal vertices is called the {\it
degree} of the tree ($n=5$ in our example). 

There is only one p.b.tree with one leaf: $\vert$ , and only one p.b.tree with two leaves: $\arbreA$.
 The main operation on p.b.trees is called {\it grafting }. The
grafting of $x$ and $y$, denoted by $x\vee y$, is the tree obtained by joining
the two roots to a new vertex:  
$$ x\vee y := \arbreAa
$$
\M

For instance: $\vert \vee \vert
= \arbreA\, , \quad \arbreA \vee \arbreA = \arbretrois\, .$
Observe that the degree of $x\vee y$ is ${\rm deg}\, x + {\rm deg}\, y +1$. 

The main point about p.b.trees is the following: the decomposition 
$$x=x^l \vee x^r$$
 exists (provided that $x\neq \vert$ ) and is {\it unique}. 
Moreover one has ${\rm deg}\, x^l < {\rm deg}\, x$ and ${\rm
deg}\, x^r < {\rm deg}\, x$.
From this property of the p.b.trees it is clear that there is a one-to-one correspondence between the
Catalan set $Y_n$, as defined in 1.1, and the set of p.b.trees of degree $n$. 

In low degree one has: 
\M

\N $Y_0 =\{\  |\  \}\ ,\  Y_1= \{\  \arbreA \  \}\ ,\ Y_2=  \{\   \arbreB
,\arbreC \  \}\ ,  \ 
Y_3=  \{\  \arbreun ,\arbredeux ,\arbretrois ,\arbrequatre ,\arbrecinq \ 
\}.$
\medskip
The union of all the sets $Y_n, n\geq 0$, is denoted $Y_{\infty}$.
\M

\N {\bf 1.3. Grove (bosquet).} In the sequel we will deal with the subsets of $Y_n$. By definition a {\it binary grove} of
degree $n$ (or simply a grove when the context is clear) is a  non-empty subset of $Y_n$. We will refer to a grove as a
disjoint union of trees.  Hence a grove is a non-empty union of p.b.trees of the same degree, such that each tree appears at
most once. We denote the set of groves of degree $n$ by $\YY_n$. The number of elements of $\YY_n$ is  $2^{c_n}-1$.
\M

\N Example: $\YY_0 = \{\vert \} ,\quad \YY_1 = \{\arbreA \} ,\quad  \YY_2 =
\{\arbreB,\, \arbreC, \, \arbreB\cup\arbreC \} .$

An important role is going to be played by the following  peculiar grove:
 $${\underline n}:= \bigcup_{y\in Y_n} y .$$
We call it the {\it total grove} of degree $n$.
 \M 
\N {\bf 1.4. Partial order structure on $Y_n$ .}
 We put  a {\it partial order} structure on the set $Y_n$ of
planar binary trees of degree $n$  as follows. We say that $x < y$ (also
denoted
$x \to y$) if the tree $y$ is obtained from $x$ by moving an edge from left to right
over a vertex, like in the basic example: 
$$\arbreB \to \arbreC.$$
 The partial order relation
is induced by this relation. More formally  the partial order on $Y_n$ is induced by the following relations:
$$\matrix{
(x\vee y) \vee z \leq x\vee(y\vee z),\cr
x < y \Longrightarrow x\vee z < y \vee z,\cr
x < y \Longrightarrow z\vee x < z\vee y .\cr
}$$

 For $n=3$
we obtain the following classical poset (pentagon):
$$\matrix{
 & &\arbreun& & \cr
 & \swarrow & &\searrow & \cr
 & & & &\arbredeux \cr
 \arbretrois& & & &\downarrow \cr
 & & & &\arbrequatre \cr
 &\searrow & &\swarrow & \cr
 & &\arbrecinq & & \cr
}$$
\centerline {figure 2}
\M

One can show that, equipped with this poset structure, $Y_n$ is a lattice, sometimes called the Tamari lattice in the
literature.
\M

\N {\bf 1.5. The Over and Under operations.} Let us introduce two new operations on planar
binary trees. For
$x\in Y_p$ and
$y\in Y_q$ the tree $x/y$ (read $x$ {\it over} $y$) in 
$Y_{p+q}$ is obtained by identifying the root of $x$ with the most left leaf of $y$.
Similarly, the  tree $x\backslash  y$ (read $x$ {\it under} $y$) in 
$Y_{p+q}$ is obtained by identifying the most right leaf of $x$ with the root of $y$.
\S

\N Examples: $\arbreC /\ \arbreA = \arbredeux$ and $ \arbreB \b \   \arbreA = \arbretrois$.

Observe that both operations are associative. The tree $\vert$ is the  neutral element for both operations and on both
sides since
$$x/\vert =\vert /x = x=
 \vert \backslash x= x\backslash \vert .$$
We will also use the following immediate property
$$x/(y\vee z) = (x/y)\vee z \quad \hbox {and}\quad  x\vee (y\b z) = (x\vee y)\b z .$$

\N {\bf 1.6. Proposition.} {\it For any p.b.trees $x$ and $y$ one has 
$$x/y\  \le\  x\backslash  y\ .$$}
\N {\it Proof.} The proof is by induction on the degree of $y$. If $y= \vert$, then both elements are equal to $x$.
If  the degree of
$y$ is strictly positive, then one can write $y=y^l\vee y^r$. Since, by induction hypothesis, $x/y^l \le x\backslash 
y^l $ and since $x/y = (x/y^l)\vee y^r$  one gets
$$x/y = (x/y^l)\vee y^r \le  (x\backslash  y^l )\vee y^r. $$
In  $(x\backslash  y^l )\vee y^r$ we move all  the vertices standing on the right leg of $x$ from left to right over the
lowest vertex to obtain $x\backslash (y^l \vee y^r)$
\BB
$$\arbreunsix$$
\B
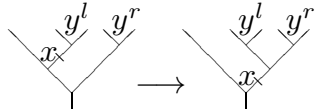
\centerline {figure 3}
\M
Since the latter tree is $x\backslash y$, we have finished the proof. \hfill $\square$
\BB

\N {\bf 2. Addition.} In this section we define a binary operation  on planar binary trees and we extend it to groves.
We call it {\it addition} or {\it sum} though it is not commutative.
\M

\N {\bf 2.1. Definition of addition.} By definition the {\it sum} of two p.b.trees $x$ and
$y$ is the following disjoint union of p.b.trees

$$ x + y := \bigcup_{x/y \le z \le x\backslash  y} z\ . $$
All the elements in the sum have the same  degree which happens to be $\deg x + \deg y$.
The associativity of the sum follows immediately from the associativity of the Over and
Under operations. 
\M

\N Example: since $\arbreA /\  \arbreA = \arbreB$ and  $\arbreA \b \ 
\arbreA =
\arbreC$, one gets 
$$\arbreA + \arbreA =  \arbreC
\cup \arbreB\, .$$

 Addition is extended to groves by distributivity on both
sides:  $$(\cup_i\, x_i)
+ (\cup_j\,y_j) := \cup_{ij}\, (x_i+y_j).$$ 
It is clear that the tree $\vert$ is the neutral element for $+$ since it is the neutral element for both the Over
operation and the Under operation.  
\M

Despite the notation $+$ that we use, the addition on groves is not commutative, for instance
$$\arbreA + \arbreB = \arbreun \ \cup \ \arbredeux \ \cup\  \arbrequatre \quad \hbox {and}
\quad
 \arbreB +  \arbreA =  \arbreun \ \cup \ \arbretrois .$$
However there is an involution, as we will see below.
\M

\N {\bf Notation.} From now on  we often denote the tree $\vert$ by 0 and the tree $\arbreA$ by 1.
\M

\N {\bf 2.2. Theorem.} {\it The sum of two groves (a fortiori the sum of two p.b.trees) is still a grove:
$$+ : \YY_n \times \YY_m \to  \YY_{n+m} .
$$}
\N {\it Proof.} Observe that it is not immediate a priori that the trees appearing in the union defining the sum  are all
different. 

This theorem is a consequence of Proposition 2.3 below.  Indeed, any grove is a subset of the
total grove. Since, in the sum of two total groves, a given tree appears at most once, the same
property is true for the sum of any two groves. Hence this sum is a grove.
\hfill
$\square$
\M

\N {\bf 2.3. Proposition.} {\it Let ${\underline n}:= \bigcup_{y\in Y_n} y$ be the total grove
of degree $n$. Then one has 
$$ {\underline n}+{\underline m}={\underline {n+m}}.$$}
\M
\N {\it Proof.} It is sufficient to prove the Proposition for $m=1$. Indeed, by induction and associativity of the
addition, it comes
$${\u n} + {\u m} = {\u n} + {\u 1} + \cdots +{\u 1} = {\u {n+1}} + {\u 1}+\cdots +{\u 1} = \cdots = {\u {n+m}} .$$

Let us show that $ {\underline n}+{\underline 1}={\underline {n+1}}$. 

We want to prove that 
$$\bigcup_{y\in Y_n}(y+1) = \bigcup_{z\in Y_{n+1}} z\ .$$
Therefore it is sufficient to show that for any element $z\in Y_{n+1}$ there exists a unique element 
$y\in Y_n$ such that  $y/1\le z \le y\b1$. We first prove the existence of $y$. We work by induction on
the degree of $z$. Let $z= z^l\vee z^r$. By induction we know that, if $z^r\neq 0$, then there exists $t$
such that $t/1 \le z^r
\le t\b 1$. One has the following relations
$$(z^l\vee t)/1 = (z^l \vee t )\vee 0 \le  z^l \vee ( t \vee 0 ) = z^l\vee (t/1) 
\le z^l\vee z^r \le z^l\vee(t\b 1)= (z^l\vee t)\b 1 .$$
Therefore one can take $y= z^l \vee t$.

If $z^r=0$, then $z=z^l\vee 0= z^l/1$. Since $z^l/1\le z^l\b 1$, $y=z^l$ is a solution.
\M
Let us show now the uniqueness of the solution. For any $y\in Y_n$ let $E(y) = \{ z\in Y_{n+1}
\mid y/1 \leq z \leq y\b 1\}$. We will show that  $\sum_{y\in Y_n} \# E(y) = c_{n+1}$.
Therefore, since any element $z\in Y_{n+1}$ belongs to some $E(y)$, it cannot belong to two of
them (since $\# Y_{n+1}= c_{n+1}$). Hence the union $\cup_{y\in Y_n} E(y)$ is a disjoint union
which covers $Y_{n+1}$.

Let $y= y^l\vee y^r$. One has
$$y/1 = ( y^l\vee y^r)/1 <  y^l\vee (y^r/1)$$
and there is no element between these two. Hence $\# E(y) = \# E(y^r) + 1$. Working by
induction, we suppose that $\sum_{y\in Y_i} \# E(y) = c_{i+1}$ for $i<n$. By using the
decomposition $Y_n=Y_{n-1}\times Y_0 \cup \cdots \cup Y_{n-i}\times Y_{i-1}
  \cup \cdots \cup Y_0\times Y_{n-1}$ we get
$$\eqalign{
\sum_{y\in Y_n} \# E(y) &= \sum_{i=1}^{n} c_{n-i} \times \big(\sum_{z\in Y_{i-1}} (\#
E(z)+1)\big)\cr &=  \sum_{i=1}^{n} c_{n-i} \times (c_i + c_{i-1}),\qquad \hbox {by induction
hypothesis},\cr &= c_n + c_{n-1}\times c_1 + \cdots + c_0\times c_n, \qquad \hbox {by
formula (1.1.1)}\cr &= c_{n+1} ,  \qquad \hbox {again by formula (1.1.1) since } c_0=1.\cr
}$$
\hfill $\square$
\M


\N {\bf 2.4. Corollary.} {\it Let $n$ and $m$ be two integers and let $z\in Y_{n+m}$. Then there
exist unique elements $x\in Y_n$ and $y\in Y_m$ such that $x/y\leq z \leq x\b y$.}\hfill
$\square$
\M

\N {\bf 2.5. Remark.} As mentioned in the introduction the elements of the Catalan set $Y_{n+1}$  are
in one-to-one correspondence with the vertices of the Stasheff polytope ${\cal K}^n$
(associahedron) of dimension $n$. One can view ${\cal K}^n$ as a cylinder ${\cal K}^{n-1}\times
I$ such that the vertices in ${\cal K}^{n-1}\times \{0\}$ (resp. ${\cal K}^{n-1}\times \{1\})$
correspond to the elements $y/1$ (resp. $y\b 1$). Then the elements between $y/1$ and $y\b
1$ are lying on the edge joining them. In particular this subset is totally ordered.
\M

\N {\bf 2.6. Involution.} Observe that for a p.b.tree {\it symmetry} around the axis passing
through the root defines an involution $\ss$ on $Y_n$ and therefore also on
$\YY_n$. For instance $\ss(\arbreA) = \arbreA$ and $\ss (\arbreB) = \arbreC$. It is clear that 
$$\eqalign{
\ss(x \vee y) &= \ss(y )\vee \ss(x),\cr
\ss(x / y) &= \ss(y )\b \ss(x),\cr
\ss(x \b y) &= \ss(y )/ \ss(x),\cr
}$$
therefore
$$\ss(x+y) = \ss(y) +\ss(x).$$

Summarizing what we have proved until now, we get the following 
\M

\N {\bf 2.7. Corollary.} {\it The set  $\YY_{\infty} := \cup_{n\geq 0}\, \YY_n$ of
groves is an involutive  graded monoid for $+$. The maps 
$$\displaylines{
\NN\longrightarrow \YY_{\infty} \longrightarrow \NN\cr
n\mapsto \u n \quad , \quad y\mapsto {\rm deg}\  y\cr
}$$
are morphisms of monoids.}
\BB

\N {\bf 3. Permutation-like notation of trees.} Though the trees help in figuring out the operations on the elements
of $Y_n$, when one wants to work explicitly it is like working with Roman numerals. So there is a need for a
more useful notation. The following permutation-like notation permits us to code linearly the elements of $Y_n$ and
so to use computer computation. It is similar to the decimal notation for integers.
\M

\N {\bf 3.1. Definition.} By definition the {\it name} of the unique element of $Y_0$ is
$0$, and the name $w(x)$ of the element $x\in Y_n$,
$n\geq 1$, is a finite sequence of strictly positive integers obtained inductively as follows. 

If $x = (y,z)\in
Y_{n-i}\times Y_{i-1}\subset Y_n$, then 
$$w(x):= (w(y), n, w(z))$$
with the convention that we do not write the zeros. If there is no possibility of confusion we simply write 
$w(x):= w(y)\ n\ w(z)$ (concatenation). Observe that, except for $0$, such a sequence is made of $n$ integers and the
integer $n$ appears once and only once.  The name of the unique element of $Y_1$ is therefore $1$, and this is in
accordance with our previous notation.
\M

\N {\bf 3.2. Bijection with the planar binary trees.} When an element in $Y_n$ corresponds to a p.b.tree $x$ and to a
name $w(x)$, we will say that $w(x)$ is the name of the tree $x$. In low dimension we get

$$\matrix{ x = &\vert   & \arbreA& \arbreB & \arbreC &
\arbreun & \arbredeux & \arbretrois & \arbrequatre & \arbrecinq \cr
w(x) = &0& 1&12 & 21 & 123 & 213 & 131 & 312 & 321 \cr
}$$

The relationship with grafting is obviously given by
 $$w(x) =\   w(x^l)\, n\,  w(x^r)\quad {\rm for}\ x=x^l\vee x^r \in Y_n.$$
Recall that symmetry around the root axis induces an involution on $Y_n$. If $a_1 \cdots a_n$ is the name of the tree
$x$, then the name of $\ss(x)$ is $a_n \cdots a_1$.
\M

\N{\bf 3.3. Weight.} From the picture of a tree  one can get its name quickly as
follows. Define its {\it
$i$th weight}  as being the degree minus 1 of the minimal subtree  which contains the leaves
number $i-1$ and $i$ (we number them from left to right starting with $0$). Then the name of
the tree is precisely the sequence of weights.
\M

\N{\bf 3.4. Test for sequences.} Given a sequence of integers, is it the name of a
p.b.treeÊ? First, check that the largest integer in the sequence is the length of
the sequence. Second,  check that the left and right subsequences (not containing the
largest integer) are names of trees. Example: 15321812 is the name of a tree, but 
15121812 is not.

 This notation has the following nice feature. If one draws the
tree metrically (leaves regularly spaced, edges at 45 degrees angle),  then the
 integers stand for the levels of the internal vertices. Check this
fact on the tree in figure 1 whose  name is $13151$.
\M

\N {\bf  3.5. Relationship with the symmetric group.} Any permutation $\ss = \ss (1) \cdots
\ss (n) $  can be made uniquely into the name of a tree by the following algorithm. Keep
the largest integer in place. On each side replace the largest integer by the length of
the side, on so on. For instance $123$ stays the same $123$, however $132$ becomes $131$,
and $23154$ becomes $13151$. The properties of the map $\SS_n \to Y_n$ so obtained and its
relationship with the polytope decompositions of the sphere as permutohedron and Stasheff
polytope respectively has been investigated in [LR1].
\BB

\N {\bf 4. The universal property.} First, we give a recursive formula for computing the sum of
groves. Second, we state the universal property of the addition which will
 enable us to construct the multiplication. 
\M

\N {\bf 4.1. Theorem [LR2].} {\it For any p.b.trees $x$ and $y$ different from $\vert$ (i.e. $0$)
the following formula holds
$$ x + y = x^l\vee (x^r + y)\, \cup \, (x+y^l)\vee y^r\, ,$$
where $x= x^l\vee x^r$ and $y=y^l\vee y^r$.}
\S
The proof follows from [LR2, theorem 5.1], where it is written in terms of associative algebras
as explained in the next section.
\M

\N {\bf 4.2. The Left and Right operations.}  From the preceding theorem it is immediately seen
 that the
sum of two trees is given as the union of two groves. In order to identify the two parts
we define two operations $\l$ and $\r$ as follows:
$$\displaylines{
x\l y := x^l\vee (x^r + y)\quad \hbox {when } x\neq 0, \cr
x\r y :=(x+y^l)\vee y^r\quad \hbox {when } y\neq 0,\cr
}$$
and 
$$ 0\l x = 0 = x\r 0.
$$
So we have
$$x+y = x\l y \, \cup \,  x\r y .$$
(Pictorially the sign $+$ splits into the two signs $\l$ and $\r$). Observe that $0\l 0 $ and $0\r
0$ are not defined though $0+0=0$, see below Remark 9.3.
\S

The operations $\l$ and $\r$ are called respectively the {\it
Left } and the {\it Right }
sum (according to the direction in which they point to). For instance one gets
$$\displaylines{
\arbreA \l \arbreA =\arbreC\ ,\quad  {\rm i.e. }\quad \ 1\l 1 = 21 ,\cr
\arbreA \r \arbreA =\arbreB\ ,\quad  {\rm i.e. }\quad \ 1\r 1 = 12 .\cr
}$$
 They are
extended to groves by distributivity with respect to disjoint union
$$ \cup_i x_i \l \cup_j y_j := \cup_{ij} (x_i \l y_j)$$
and similarly for $\r$.

Since the Left (resp. Right)  sum of two groves is a subset of the  sum, the Left (resp. Right)  sum of
groves is a grove.

It is clear that the relationship with the involution is as follows
 $$\ss(x\l y) = \ss (y)\r\ss (x)\quad {\rm and}\quad   \ss(x\r y) = \ss (y)\l\ss (x).\eqno (4.2.1)$$

Recall that $x/y = (x/y^l)\vee y^r$ and $x^l\vee (x^r\b y)= x\b y$. One can show that in terms of the
operations $/$ and $\backslash$ we get
$$x\l y =  \bigcup_{x^l\vee (x^r/ y)\le z \le x\backslash y} z \quad
\hbox {and }\quad x\r y =
  \bigcup_{x/  y \le z \le (x\b y^l)\vee y^r} z
\eqno (4.2.2)$$

In the Appendix the addition table is written such that the first line gives $x\r
y$ and the second (union third if any) line gives $x\l y$.
\M

\N {\bf 4.3. Tricks for computation.} For some special trees the computation of the
Left or Right sum is easy. First recall that  $-\vee - = - {\tt max} -$, where {\tt max} stands
 for the largest integer (i.e the length) of the word to which it pertains. For instance
$1\vee 1 =131$ and $1\vee 0= 12$. Then it is easy to check that
$$\displaylines{
(- {\tt max}) \l - = - {\tt max }-  \cr 
-  \r ({\tt max} -) = - {\tt max }-  \cr 
}$$
For instance:
$$213 \l 12 = 21512,\quad 1412\r 54131 = 141294131.$$

\N {\bf 4.4. Proposition.} {\it The Left and Right sum on groves satisfy the following relations
$$\left\{\eqalign{
(x\l y) \l z = x\l (y + z),\cr
(x\r y) \l z = x\r (y \l z),\cr
(x + y) \r z = x\r (y \r z).\cr
}\right.$$
and $0\r x = x = x\l 0$ for $x\neq 0$.}
\M
The proof follows from [L2] as we will show in section 5.
\M
\N {\bf 4.5. Theorem.} {\it For any p.b.tree $x$ of degree $n$ there is a
unique way of writing it as a composition of $n$ copies of 1 with the Left
and the Right sum. 
It is called the ``universal expression" of $x$ and
denoted $w_x(1)$.

The universal expression is unique modulo the relations of 4.4.}
\M

The proof follows from [L2] as we will show in section 5.
Observe that this theorem gives still another combinatorial description of the Catalan sets (not
in the list of [St]).

The inductive algorithm to construct $w_x(1)$ is given by $w_0(1)= 0$ and:
$$w_x(1):= w_{x^l}(1)\r 1 \l w_{x^r}(1).\eqno (4.5.1)$$

\N Examples:
$$\eqalign{
12 &= 1\r 1\cr
21 &= 1 \l 1\cr
123 &= 12 \r 1 = (1\r1)\r 1\cr
213 &= 21 \l 1 = (1\l 1)\r 1\cr
131 &= 1\r 1 \l 1 \cr
312 &= 1\l 12 = 1\l (1\r 1)\cr
321 &= 1\l 21 =1\l (1\l 1)\cr
131492141&= ((1\r 1\l 1) \r 1) \r 1\l ((1\l 1) \r 1\l 1).\cr
}$$
\M
\N {\bf 4.6. Remark.} Let $a_1 a_2 \cdots a_n$ be the name of the tree $x$. We say that there is
an {\it ascent}  (resp. a {\it descent}) at $i$ if $a_i< a_{i+1}$ (resp.  $a_i>a_{i+1}$). One can show that, in
the universal expression of $x$, the signs are $\r$ when there is an ascent and $\l$ when there
is a descent.
\B
\N {\bf 5. Polynomial algebra with tree exponents.} In this section we show how the results of
section 4 are just translation of results contained in [LR2], [LR1] and [L2].

Let $K$ be a commutative ring. It is well-known that the polynomial ring $K[X]$
has a linear basis indexed by $\NN$:
$\{X^n\}_{n\in \NN}$. Multiplication of monomials in $K[X]$ corresponds to
addition in $\NN$: $X^n X^m = X^{n+m}$. Observe that if we were using the group algebra notation,
then the polynomial algebra would be denoted $K[\NN]$.

Similarly let us introduce $K[Y_{\infty}]$ the vector space with basis $X^y,
y\in Y_n$ for all $n\ge 0$. We put a product on it by the formula:
$$X^y * X^z := X^{y+z}.$$
Since $y+z$ need not be a tree but a grove, we define
 $X^{y\cup y'}$ as:
$$X^{y \cup y'}:= X^{y}+X^{y'},$$
so that the previous formula has meaning. Under  this
formula, any grove determines an element in $K[Y_{\infty}]$. This associative algebra has
already been encountered in the framework of dendriform algebras as we now explain.
\M
\N {\bf 5.1. Dendriform algebras [L1, L2].} By definition a {\it dendriform algebra} (also called {\it dendriform dialgebra})
is a
$K$-vector space  $A$ equipped with two binary operations
$$
\eqalign{
&\prec\ : A\otimes A\to A,\cr
&\succ\ : A\otimes A\to A,\cr}
$$
which satisfy the following axioms:
$$
\eqalign{
\hbox{(i)}\quad &(a\prec b)\prec c = a\prec(b*  c)  ,\cr
\hbox{(ii)}\quad& (a\succ b)\prec c = a\succ (b\prec c)  ,\cr
\hbox{(iii)}\quad& (a* b)\succ c = a\succ (b\succ c)  .\cr}
$$
for any elements $a, b$ and $c$ in $E$, with the notation
$$a* b:= a\prec b + a \succ b .$$
Adding the three relations show that the operation $* $ is associative.
\M

\N {\bf 5.2. A dendriform algebra associated to p.b.trees.} Let  $K[Y_{\infty}']$ be the vector
space generated by the elements
$X^x, $ for $x\in Y_n\  (n\geq 1)$ (we exclude $0\in Y_0$ for a while). We define two operations
on  $K[Y_{\infty}']$ by the formulas
$$\eqalign{
X^x \prec X^y := X^{x^l\vee(x^r+y)} ,\cr
X^x \succ X^y := X^{(x+y^l)\vee y^r)} ,\cr
}$$
with the convention that, for any two trees $y$ and $y'$, one has 
$$X^{y \cup y'}:= X^{y}+X^{y'}.$$
These operations are extended to $K[Y_{\infty}']$ by distributivity:
$$ (a+b)\prec c = a\prec c + b\prec c \ {\rm and }\ a\prec (b+c) = a\prec b + a\prec c $$
and similarly for $\succ$.
\M

\N {\bf 5.3. Theorem (Universal property)} [L2, Proposition 5.7].  {\it The vector space $K[Y_{\infty}']$ 
equipped with the two operations $\prec$ and $\succ$ as above is a dendriform algebra.
Moreover it is the free dendriform algebra on one generator (namely $X$).}\hfill $\square$
\M

\N {\bf 5.4. Proof of the results of section 4.} The associative product on $K[Y_{\infty}']$ is
given by 
$$X^x *X^y = X^x\prec  X^y + X^x\succ  X^y = X^{x\l y} + X^{x \r y} = X^{x+y}.$$
It is a nonunital associative algebra. We add a unit  $1= X^0$ to it, so that $K[Y_{\infty}]= K[Y_{\infty}']
\oplus K\cdot1$
becomes an augmented unital associative algebra. The operations  $\prec$ and $\succ$ can
 be partially extended to $K[Y_{\infty}]$ by  $X^0\succ X^x = X^x = X^x\prec X^0$ for $x\neq 0$.
\M

In [LR2, theorem 5.1] it is proved that this algebra structure on $K[Y_{\infty}]$ is the same as the one given
by the rule
$$X^x *X^y = \sum _{x/y\leq z \leq x\b y} X^z.$$
Therefore the results of section 4 are just translations of the results of [LR2] and the results of
[LR1] cited above.
\M
\N {\bf 5.5. Remarks about notation.} In [LR1, LR2, L1, L2] the linear generators are denoted $x$
instead of $X^x$. We adopt this new notation here to avoid confusion with the operations in
$\YY_{\infty}$.

In [L1] and [L2] the symbols $\l$ and $\r$ are used to denote operations of an associative
dialgebra which is the Koszul dual structure of dendriform algebra. We have given them a
completely different meaning here.
\B

\N {\bf 6. Multiplication.} Since the polynomial algebra is the free associative algebra on one
generator, one can define the composition of polynomials. It turns out that composite of
monomials is still a monomial. It is related to multiplication of integers by: $(X^n)^m = X^{nm}$.

Similarly, since the associative algebra $K[Y_{\infty}']$ is free on one generator when
considered as a dendriform algebra, one can perform composition of polynomials with tree
exponents. Though the composite of monomials, that is $(X^x)^y$ where $x$ and $y$ are
p.b.trees, is not a monomial, it turns out that it is $X$ to the power of some grove. Hence one
can define the multiplication of p.b.trees as being this grove and then extend this
multiplication to any groves.
\M
\N {\bf 6.1. Definition.} Let   $x$ and $y$ be planar binary trees. By definition the {\it product}
$x\times y$ is
$$ x \times y := w_x(y).$$
where $ w_x(1)$ is the universal expression of $x$ (cf. Theorem 4.5). In other words we replace
all the copies of $1$ by copies of $y$ in this universal expression.

Observe that the above definition of the product has a meaning even when $y$ is a grove since
Right sum and Left sum of groves are well-defined.  We extend the multiplication to $x$ being a
grove by distributivity on the left with respect to disjoint union:
$$(x \cup x')\times y = w_x(y)\cup w_{x'}(y).$$
So we have defined the product of two groves. This product is a grove since it is obtained by the
operations $\l$ and $\r$. It is clear that the degree of the product is the product of the degrees, so we have defined
a map:
$$ \times : \YY_n \times \YY_m \to  \YY_{nm} .
$$
Observe that the product is not commutative.
\M

\N {\bf Examples.} Since $21=1\l 1$, one has $21\times x = x\l x$, so, for
instance, $21 \times 12 = 12\l 12 = 1412$. On the other hand
$12\times 21 = 21 \r 21 = 2141$. 

Since $131 = 1\r 1 \l
1$,  one has $131\times x = x\r x\l x$, so, for instance, $131 \times 21 = 12\r
12\l 12 = 1412\l 12= (1412+1)\vee 0 = \allowbreak (1412\l 1)6 \cup 141256=
((1\vee (12+1))6 \cup 141256= 151316 \cup 151236\cup 141256$.
\M

 \N {\bf 6.2. Proposition.} {\it The multiplication $\times$ on groves is
distributive with respect to the Left sum, the Right sum and the sum on the left, but not on the right.}
\M
\N {\it Proof.} The formula $w_{x+x'}(1) = w_x(1)+ w_{x'}(1)$ follows from 
$$\eqalign{
w_{x\l x'}(1) &=  w_x(1)\l w_{x'}(1),\ {\rm and }\cr
w_{x\r x'}(1) &=  w_x(1)\r w_{x'}(1).\cr
}$$

These last two formulas follow inductively from the properties of the function $w$, namely
(4.5.1) and $w_1(1)=1$.\hfill $\square$
\M

 \N {\bf 6.3. Proposition.} {\it
The multiplication of groves is associative with neutral element on both sides the tree
$\arbreA = 1$.} 
\M
\N {\it Proof.}  Interpreted in terms of dendriform algebra, the multiplication of p.b.trees is
composition of monomials. Indeed, since $K[Y_{\infty}']$ is the free dendriform algebra on the generator $X$, there exists a unique morphism of
dendriform algebras $W_x : K[Y_{\infty}'] \to K[Y_{\infty}']$ sending $X$ to $X^x$. The image of $X^y$ by $W_x$ is precisely
$X^{w_x(y)}$. Since composition of dendriform algebra morphisms is associative, the multiplication of p.b.trees is
associative.

Since $w_1(1)=1$, we get $w_1(y)=y$ and so $1\times y = y$. On the other hand $x\times1 = x$
is a tautology. \hfill $\square$
\M

\N {\bf 6.4. Theorem.} {\it With the notation ${\u  n} = \cup_{y\in Y_n}\ 
y$, one has 
$${\u  n} \times {\u  m} = {\u  {nm}}.$$}

\N {\it Proof.} Since the multiplication is distributive on the left with respect to the addition, we
get 
$$\eqalign{
{\u  n} \times {\u  m} &= ({\u  1}+ \cdots + {\u  1})\times {\u  m} = {\u  1}\times {\u  m}
+ \cdots + {\u  1}\times {\u  m}\cr
&= {\u  m} + \cdots + {\u  m}=  {\u  {nm}},\cr
}$$
since ${\u  1}\times {\u  m} = {\u  m}$ by Proposition 6.3.
 { }\hfill $\square$
\M

\N {\bf 6.5. Proposition (recursive property).} {\it Let $x=x^l\vee x^r$ be a p.b.tree and let $y$ be
a grove. The multiplication is given recursively by the formulas
$$x\times y = (x^l \times y) \r y \l (x^r \times y)$$
and $0\times y = 0$. }
\M
\N {\it Proof.} The universal expression satisfies
$$w_{x^l\vee x^r}(1) = w_{x^l}(1) \r 1 \l w_{x^r}(1) $$
because the name of $x^l\vee x^r$ is $w(x^l)\ n \ w(x^r)$ (cf. 4.3).  Hence we get
$$\eqalign{
w_{x^l\vee x^r}(y) &= w_{x^l}(y) \r y \l w_{x^r}(y)\cr
&= (x^l\times y) \r y \l (x^r\times y).\cr
}$$
\hfill $\square$
\S
\N {\bf Exercise.} Starting with the formula of Proposition 6.5 as a definition for $\times $, prove associativity. Hint: use
the formula 
$$
a\vdash b \dashv c = (a + b^l) \vee (b^r + c).
$$
\S
\N {\bf 6.6. Proposition (involution).} {\it For any groves $x$ and $y$ one has
$$\ss(x\times y) = \ss (x) \times \ss (y).$$}
\M
\N {\it Proof.} It obviously suffices to prove the formula when $x$ is a p.b.tree. We work by
induction on the degree of $x$. The formula is a tautology for $x=1$. By Proposition 6.5 we get
$$\eqalign{
\ss(x\times y) &= \ss (x^l \times y \vdash y \dashv x^r\times y)\cr
&=\ss(x^r\times y)\r \ss(y)\l \ss(x^l\times y) \cr
&= \big(\ss(x^r)\times \ss (y)\big)\r \ss(y)\l \big(\ss(x^l)\times \ss (y)\big)\cr
&= \big( \ss(x^r)\vee \ss(x^l)\big)\times \ss (y)\cr
&= \ss(x)\times \ss (y).\cr
}$$
\hfill $\square$

\N {\bf 6.7. Summary.} On the set of groves $\YY_{\infty} =
\bigcup_{n\ge 0}\YY_n$ there are defined operations + and $\times$ such that

$\bullet$ the addition + is associative, distributive both sides with respect
to $\cup$, with neutral element $0=\vert$ , but is not commutative,

$\bullet$ the multiplication $\times$ is associative, distributive on the left
with respect to the sum + and to the disjoint union $\cup$ (but not right distributive), with neutral
element (both sides) $1=\arbreA$, but is not commutative,

$\bullet$ the involution $\ss$ on $\YY_{\infty}$ satisfies $\ss(x+y) = \ss (y) + \ss (x)$ and 
$\ss(x\times y) = \ss (x) \times \ss (y)$,

$\bullet$  the maps $\NN \to \YY , n\mapsto {\underline n}=\bigcup_{y\in Y_n}y$ and 
 $\deg : \YY \to \NN$
(degree) are compatible with $+$ and $\times$. The composite is the identity of
$\NN$.
\B

\N {\bf 6.8. Questions and problems.} By definition a grove is {\it prime}
if it is not the product of two groves different from 1. Obviously
any tree (resp. grove) whose degree is prime is a prime tree (resp. grove). However
 there are also prime trees and
groves of nonprime degree, for instance $1234$. The tree $1241$ is not prime since
 $1241 = 12\times 21$.

It would be interesting to study the factorization of a grove by another one. In particular it seems that, when a grove is a
product of prime groves, then the ordered sequence of factors is unique.
\B

\N {\bf 7. Elementary combinatorial applications.} Here are simple combinatorial applications of the formulas ${\u n} + {\u
m} = \u {n+m}$ and ${\u n} \times {\u m} = \u {nm}$ respectively.
\M
\N {\bf 7.1. Proposition.} {\it  For any pair of p.b. trees $x\in Y_n$ and 
$y\in Y_m$ let $c_{x,y}$ be the  number of trees in the grove $x+y$. Then
 the
following combinatorial formula holds:
$$ \sum_{x\in Y_n,\, y\in Y_m}c_{x,y} = c_{n+m}.
$$}
\hfill $\square$
\M

We see for instance that $c_{12,12} = 3,\ c_{12,21} = 2,\ c_{21,12} = 6,\ c_{21,21} = 4, 
$ and so $3+2+6+3 = 14 = c_4$ as expected.
\M

\N {\bf 7.2. Proposition.} {\it  For any   p.b. trees $x\in Y_n$ and 
any integer $m$  let $d_{x,{\u m}}$ be  the number of trees in the grove  $x\times {\u m}$. Then
 the
following combinatorial formula holds:
$$ \sum_{x\in Y_n}d_{x,{\u m}} = c_{nm}.
$$}
\hfill $\square$

We see, for instance that
\S

\N  $d_{12,{\u 2}} = 7, d_{12,{\u 2}} = 7,  
$ and so $7+7 = 14 = c_4$ as expected, 
\S

\N  $d_{123,{\u 2}} = 22,\ d_{213,{\u 2}} = 33,\  d_{131,{\u 2}} =20,\ d_{312,{\u 2}} = 33,\ d_{321,{\u 2}} =22,
$ and so $22+33+20+33+22 = 132 = c_6$ as expected.
\BB
\vfill
\eject

\N {\bf Appendix to part I: tables for planar binary trees.}
\B

\N {\bf I.A.1. Addition table.} Recall that $0$ is the neutral element for $+$ , so
$$0+x = x = x+ 0.$$
In the following table we omit the $\cup$ sign, $x/y$ is the first element of the first
line, $x\backslash y $ is the last element of the last line. The first line is $x\r y$ and the
second line (union third if any) is $x\l y$. Recall that $\ss (x+y) = \ss (y) + \ss (x)$, so a sum
like $21 + 123$ is easily obtainable from this table.
\B

{\vbox{
\hrule
\halign{&\vrule#&\strut \hfil# \cr
$\ x + y$	&&\ 1 		&&\ 12 			&&\ 21		&& \cr
\noalign{\hrule}
\ 1  	&&\ 12 		&&\ 123\  213 			&&\ 131		&& \cr 
	&&\ 21     	&&\ 312		&&\ 321		&& \cr
\noalign{\hrule}
\ 12 	&&\ 123	    	&&\  1234\ 1314		&&\ 1241		&& \cr
	&&\ 131	&&\ 1412			&&\ 1421 	 	&&\cr
\noalign{\hrule}
\ 21 	&&\ 213 &&\  2134\ 3124\ 3214 &&\ 2141 	&& \cr 
	&&\  312	\ 321		&&\ 4123\  4213\  4312		&&\  4131\ 4321\ 		&& \cr
\noalign{\hrule}
\ 123	&&\ 1234		&&\  12345\ 12415		&&\ 12351		&& \cr
	&&\ 1241 		&&\ 12512		&&\ 12521		&& \cr
\noalign{\hrule}
\ 213		&&\ 2134	&&\  21345\ 21415		&&\ 21351			&& \cr
&&\ 2141			&&\ 21512		&&\ 21521	&& \cr
\noalign{\hrule}
\ 131	&&\ 1314		&&\  13145\ 14125\ 14215	&&\ 13151		&& \cr
	&&\  1412\ 1421 	&&\  15123\  15213\ 15312	&&\  15131\ 15321 	&& \cr
\noalign{\hrule}
\ 312	&&\ 3124		&&\  31245\ 41235\  41315 &&\ 31251 	&& \cr
	&&\  4123\ 4131	&&\ 51234\  51314\ 51412	&&\ 51241\ 51421 && \cr
\noalign{\hrule}
\ 321 &&\ 3214		&&\ 	32145\ 42135\ 43125\ 43215 &&\ 32151		&& \cr
	&&\ 4213\  4312\ 4321\ &&\ 52134\  53214\  53124\quad &&\  54131\ 
52141\ 54321\ &&\cr
 	&&		&&\qquad   54123\  54213\ 54312\  &&		&& \cr
\noalign{\hrule}
  }}
\BB
\N Example of computation by using the recursive formulas:
$$\eqalign{
131\vdash 12 &= 131 \vdash (1\vee 0) = (131 + 1) \vee 0\cr
& =( 131\l 1 \ \cup \ 131 \r 1 )\vee 0 = ((1\vee 1)\l 1 \ \cup \ 1314)\vee 0 \cr
& =((1\vee (1+1) \ \cup\ 1314)\vee 0 = ( 1412\  \cup \ 1421\ \cup 1314  )\vee 0 \cr
& =  14125\  \cup \ 14215\ \cup\  13145.  \cr
}$$
\vfill\eject

\N {\bf I.A.2. Mutiplication table.} Recall that 1 is the neutral element for $\times$, so 
$$1 \times x =x =x\times 1.$$
Recall that the recursive formula is $x\times y = (x^l \times y) \r y \l (x^r \times y)$ and
$ 0\times x = 0, \ 0\r x = x = x\l 0$.
\B
\petit{
{\vbox{
\hrule
\halign{&\vrule#&\strut \hfil# \cr
$\ x \times y$	&&\ 12 		&&\ 21		&&\ 12 \ 21		&&\ 123 &&\ 213
&&\ 131 &&\ 312 &&\ 321 && \cr
\noalign{\hrule}
\ 12  &&\ 1234	&&\ 2141 &&\ 1234\ 1314\ 1241 &&\ 123456\ &&\ 213516 &&\ 131461 &&\ 312612 &&\ 321621 && \cr 
\  	  &&\ 1314 &&\ 		   &&\ 2134\ 3124\ 3214	&&\  124156 &&\ 215216 &&\   141261  &&\        &&\        &&
\cr
\  	  &&\ 		     &&\ 		   &&\ 	2141	                   &&\  125126 &&\              &&\   142161     &&\        &&\        && \cr
\noalign{\hrule}
\ 21  &&\ 1412	&&\ 4131 &&\ 1412\ 4123\ 4213 &&\ 126123 &&\ 216213 &&\ 161241 &&\ 612512 &&\ 621521 &&\cr
\  	  &&\ 		   &&\ 4321	&&\ 4312\ 1421\ 4131	&&\              &&\              &&\ 162141 &&\ 615312   &&\  651421 &&
\cr
\  	  &&\ 		   &&\ 		   &&\ 4321		                      &&\              &&\             &&\   164131 &&\   &&\  654321      && \cr
\noalign{\hrule}
  }}

\BB

{\vbox{
\hrule
\halign{&\vrule#&\strut \hfil# \cr
$\ x \times y$	&&\ 12 		 &&\ 21		    &&\ 12  $\cup$   21 		&& \cr
\noalign{\hrule}
\ 123 &&\ 123456	&&\ 214161		&&\ 123516\ 123456\ 131516\ 131456\ 125126\ 125216		&& \cr
\  	  &&\ 131456	&&\ 		      &&\ 124156\ 213516\ 213456\ 312516\ 312456\ 321516  && \cr
\  	  &&\ 123516  &&\ 		      &&\ 321456\ 215126\ 215216\ 214156\ 123461\ 131461	&& \cr
\  	  &&\  131516  &&\ 		      &&\ 124161	\ 213461\ 312461\ 321461\ 214161	           && \cr
\noalign{\hrule}
\ 213 &&\ 141256		       &&\ 413161		&&\ 151236\ 151316\ 141256\ 152136\ 153126\ 153216		&& \cr
\  	  &&\ 151236		       &&\ 432161		&&\ 142156\ 512346\ 512416\ 412356\ 521346\ 521416 	&& \cr
\  	  &&\ 151316		       &&\ 			     &&\ 421356\ 531246\ 541236\ 541316\ 431256\ 513146  && \cr
\  	  &&         	       &&\ 			     &&\ 514126\ 514216\ 413156\ 532146\ 542136\ 543126  && \cr
\  	  &&         	       &&\ 			     &&\ 543216\ 432156\ 141261\ 142161\ 412361\ 421361  && \cr
\  	  &&         	       &&\ 			     &&\ 431261\ 413161\ 432161                          && \cr
\noalign{\hrule}
\ 131 &&\ 123612		       &&\ 216131		&&\ 123612\ 123621\ 131612\ 131621\ 126123\ 126213		&& \cr
\  	  &&\ 131612		       &&\ 216321		&&\ 126312\ 126131\ 126321\ 213612\ 213621\ 312612	 && \cr
\  	  &&        		       &&        		&&\ 312621\ 321612\ 321621\ 216123\ 216213\ 216312	 && \cr
\  	  &&         	       &&        		&&\ 216131\ 216321                                	 && \cr
\noalign{\hrule}
  }}
}
\BB
\normal

The products $312 \times y$ and $321\times y$ for $y= 12$ or $21$ or $12\, \cup\, 21$ can be obtained from this
table by the formula $\ss (x\times y) = \ss (x) \times \ss (y)$.
\M

 \N {\bf I.A.3. Exercise.} The following grove  ($\cup$
is  omitted): 
$$\matrix{
141294131&  141291241&  141292141 \cr
142194131&  142191241&  142192141\cr 
131494131&  131491241&  131492141\cr
}$$
is a product of two groves, which ones ?
\BB

\centerline {\bf II. Arithmetic of planar trees}
\B

In this second part we extend the addition and the multiplication on planar binary trees to all planar trees. The binary
case becomes a quotient of this case. The formulas are slightly more complicated since there is now one more tree in
degree 2 which accounts to one more operation that we call the Middle sum. As before this theory is governed by some
free object on one generator in a category of algebras, they are the {\it dendriform trialgebras} introduced in [LR3].

Instead of using a poset structure we begin right away with the recursive definition of the sum. We skip the arguments
which are similar to the arguments in the first part.
\B

\N {\bf 8. The set of planar  trees.} 
\M
\N {\bf 8.1. Super Catalan sets.} Let $T_0$ be a set with one element. The sets $T_n$ for
$n\geq 1$ are defined inductively by the formula
$$T_n:= \bigcup_{i_0+ \cdots + i_k = n-k} T_{i_0}\times \cdots \times T_{i_k}$$
where the disjoint union is extended to all possibilities with $k\geq 1$ and $i_j \geq  0$ for all $j$.

If we denote by $a$ the unique element of $T_0$, then an element of $T_n$ can be described as a partial parenthesizing
of the word $aa\cdots  a$ of length $n+1$. Let $x^{(0)}\in T_{i_0}, \cdots , x^{(k)}\in T_{i_k}$. The element
$(x^{(0)},\cdots , x^{(k)})\in  T_{i_0}\times\cdots \times T_{i_k}$ viewed as an element in $T_n$ is denoted
$x^{(0)}\vee\cdots \vee x^{(k)} \in T_n$. In terms of parenthesizing it simply consists in concatenating the 
words and putting parenthesizes at both ends.

There are many other combinatorial descriptions of the sets $T_n$. We will
use two of them as described below, one classical: the planar  trees, and one less classical: the
permutation-like notation. Others include the cells of the 
$(n-1)$-dimensional Stasheff polytope.

Let $C_n$ be the number of elements of $T_n$. It comes immediately: 

$C_0$ =1 and 
$$C_n =\sum_{i_0+ \cdots + i_k = n-k} C_{i_0}  \cdots C_{i_k}.$$
 It is  classically called the Super Catalan
number, so $T_n$ is called the {\it Super Catalan set} (see below for the generating series). 
\M

\N {\bf 8.2. Planar trees.} A {\it planar tree} is an oriented planar
graph drawn in the plane with $n+1$ leaves and one root, such that each internal
vertex has at least two leaves and one (and only one) root. We consider these trees up to isotopy.
Here is an example:
\M

$$\trree$$

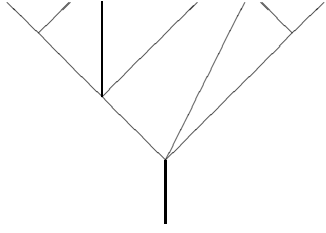
\centerline {figure 4}
\M
If the number of leaves is $n+1$, then we say that the tree is of degree $n$ (in our example $n=6$ ). 

There is only one planar tree with one leaf: $\vert$ , and only one planar tree with two leaves: $\arbreA$.
 The main operation on planar trees is called {\it grafting}. The
grafting of $x^{(0)}, \cdots , x^{(k)}$ denoted by $x^{(0)}\vee  \cdots \vee x^{(k)}$ is the tree obtained by joining
the roots to a new vertex and creating a new root:  
$$x^{(0)}\vee  \cdots \vee x^{(k)} := \qquad \arbreAaa
$$
\M

For instance: $\vert \vee \vert \vee \vert
= \arbredeuxdeux$.
Observe that the degree of $x^{(0)}\vee  \cdots \vee x^{(k)}$ is $\sum_i ({\rm deg}\,x^{(i)} + 1) -1$. 

The main point about planar trees is the following: given a planar tree $x$ the decomposition 
$$x=x^{(0)}\vee  \cdots \vee x^{(k)}$$
 exists  and is {\it unique}. 
Moreover, when the degree of $x$ is strictly positive,  one has $k\geq 1$ and ${\rm deg}\, x^{(i)} < {\rm deg}\, x$  for
any $i$. From this property of the planar trees it is clear that there is a one-to-one correspondence between the
Super  Catalan set $T_n$, as defined in 8.1, and the set of planar trees of degree $n$. So we identify them.

In low degree one has: 
\M

$$T_0 =\{\  |\  \}\ ,\  T_1= \{\  \arbreA \  \}\ ,\ T_2=  \{\   \arbreB
,\arbreC \  , \arbredeuxdeux\ \}$$
$$T_3=  \{\  \arbreun ,\cdots \ , \arbreuntroistrois \ , \cdots , \arbretroistroistrois\   \}.$$
The union of all the sets $T_n$ for $n\geq 0$ is denoted $T_{\infty}$.
\M
Each set $T_n$ is itself graded by the number of internal vertices. Denoting by  $T_{n,i}$ the set of planar trees with
$n+1$ leaves and $n+1-i$ internal vertices, we get a decomposition
$$ T_n = T_{n,1} \cup \cdots \cup T_{n,n}.$$
It is clear that $T_{n,1} = Y_n$ (planar binary trees), and that $T_{n,n}$ is made of only one element. In the Stasheff
polytope interpretation of the Super Catalan sets the elements of $T_{n,i}$ correspond to the $(i-1)$-cells of the
Stasheff polytope of dimension $(n-1)$. Let
$a_{n,i}$ be the number of elements of
$T_{n,i}$. We define a generating series by
$$ f(x,t) := 1+\sum_{n\geq 1} (\sum_{i=1}^{i=n}a_{n,i}t^{i-1}) x^n.$$
One can show, either by direct inspection like in 1.1, or by using Koszul duality of operads (cf. [LR3, LR4]), that 
$$ f(x,t) =   {1+tx -\sqrt { 1-2(2+t)x+t^2x^2 }\over 2(1+t)x}\ .$$
\S

\N {\bf 8.3. Groves (bosquets).} In the sequel we will deal with the subsets of $T_n$. By definition a {\it grove} of
degree $n$ is a  non-empty subset of $T_n$. We will refer to a grove as a disjoint union of trees.  Hence a grove is a
non-empty union of planar trees of the same degree, such that each tree appears at most once. We denote the set of
groves of degree $n$ by $\TT_n$. The number of elements of $\TT_n$ is  $2^{C_n}-1$. The set of groves made of the trees
in $T_{n,i}$ is denoted $\TT_{n,i}$.
\M

An important role is going to be played by the following  peculiar grove:
 $${\underline n}:= \bigcup_{x\in T_n} x .$$
We call it the {\it total grove} of degree $n$.
\BB

\N {\bf 9. Addition.} In this section we define a binary operation  on planar trees and we extend it to groves.
We call it {\it addition} or {\it sum} though it is not commutative.
\M

\N {\bf 9.1. Definition of addition.} By definition the {\it sum} of two planar trees $x$ and
$y$ is given recursively by the formula:
$$\eqalign{ x+ y :=&\quad   x^{(0)}\vee \cdots \vee (x^{(k)}+ y) \cr
&\cup \ x^{(0)}\vee \cdots \vee (x^{(k)}+y^{(0)})\vee \cdots \vee y^{(\ell )} \cr
&\cup \  (x+y^{(0)})\vee \cdots \vee y^{(\ell )}.\cr
}$$
where $x = x^{(0)}\vee \cdots \vee x^{(k )}\in T_p$ and $y = y^{(0)}\vee \cdots \vee y^{(\ell )}\in
T_q$ (for $p\neq 0 \neq q$).
Moreover  $0= \vert \in T_0$ is a neutral element for  $+$. 

One observes that the sum of two trees is a (disjoint)  union of
trees. Their degree is the sum of the degrees of the starting trees.

\M

\N Example:
$$\arbreA + \arbreA =  \arbreC
\cup \arbreB\cup \arbredeuxdeux \, .$$

The addition is extended to the union of trees  by distributivity on both
sides:  $$(\cup_i\, x_i)
+ (\cup_j\,y_j) := \cup_{ij}\, (x_i+y_j).$$ 
\N {\bf Notation.} From now on  we often denote the tree $\vert$ by 0 and the tree $\arbreA$ by 1.

Despite the notation $+$ that we use, the addition on groves is not commutative,
however there is an involution, as we will see below.
\M

\N {\bf 9.2. The three operations $\l\ , \ \r \ , \ \m $.}
 From the definition of the sum it is clear that it comes as a disjoint union of three
different pieces. Let us define the following three operations (for $x\neq 0$ and $y\neq 0$)
$$\eqalign{
x\l y &:=  x^{(0)}\vee \cdots \vee (x^{(k)}+ y)\ , \cr
x\r y&:=  (x+y^{(0)})\vee \cdots \vee y^{(\ell )}\ ,\cr
x\m y &:=  x^{(0)}\vee \cdots \vee (x^{(k)}+y^{(0)})\vee \cdots \vee y^{(\ell )}\ .\cr
}$$
We call $\l$ the {\it Left sum, $\r$ the Right sum,} and  {\it $\m$ the Middle sum}. They are well-defined for
planar trees and extended to groves by distributivity. By construction one has 
$$x+y := x\l y\ \cup\ x\r y \ \cup\ x\m y \ .$$

We extend these definitions to $x$ or $y$ being $0$ in the following cases: for $x\neq 0$

\qquad $\Big\lbrace \matrix {x\l 0 = x\cr
                    0\l x = 0\cr}$\qquad 
$\Big\lbrace \matrix {x\r 0 = 0\cr
                    0\r x = x\cr}$\qquad 
$\Big\lbrace \matrix {x\m 0 = 0\cr
                    0\m x = 0\ .\cr}$

For $x=0$ these elements are not defined, however $0+0=0$.
\M

\N {\bf 9.3. Remark.} In order to give a meaning to $0\l 0, 0\r 0 $ and   $0\m 0$, we could use the following trick.
Think of a planar tree as being define by its internal vertices and the relationship between them. So the tree $\vert$
is in fact the empty set $\emptyset$ since it has no vertex. Then one can put, without contradiction, 
 $\emptyset \l \emptyset =\emptyset \r \emptyset  =\emptyset \m \emptyset = \emptyset $. We still have $\emptyset +
\emptyset =
\emptyset$ since the union of the empty set with itself is still the empty set.
\M

\N {\bf 9.4. Proposition.} {\it The sum of union of trees is
associative. The neutral element if $0$.}
\M
\N {\it Proof.} 
In the next section we will show that this associativity property follows from [LR3, LR4] (cf. 11.4). \hfill $\square$
\M

\N {\bf 9.5. Theorem.} {\it  The  Left sum, the Middle sum,  the Right sum and the  sum  of two groves (a fortiori of
two planar trees) is still a grove.}
\M
\N {\it Proof.} 
This theorem is a consequence of Proposition 9.6 below.  Indeed, any grove is a subset of the
total grove. Since, in the sum of two total groves, a given tree appears at most once, the same
property is true for the sum of any two groves. Hence this sum is a grove.
\hfill
$\square$
\M

\N {\bf 9.6. Proposition.} {\it Let ${\underline n}:= \bigcup_{x\in T_n} x$ be the total grove
of degree $n$. Then one has 
$$ {\underline n}+{\underline m}={\underline {n+m}}.$$}
\M
\N {\it Proof.} It is sufficient to prove the Proposition for $m=1$. 
Let us show that $ {\underline n}+{\underline 1}={\underline {n+1}}$. 

We want to prove that 
$$\bigcup_{t\in T_n}(t+1) = \bigcup_{s\in T_{n+1}} s\ .$$
Therefore it is sufficient to show that for any element $s\in T_{n+1}$ there exists a unique element 
$t\in T_n$ such that  $s \in t+1$. 

First we show that, for a given tree $t$,  the trees appearing in $t+1$ are obtained from $t$ by adding a new leaf
starting from the right side  either from a vertex (lying on this right side) or from the middle of an edge (including the
root and the leaf). By definition $t+1 = t\l 1\ \cup \ t\m 1 \ \cup \ t\r 1$. The third component  $t\r1$  is 
$(t+0)\vee 0 = t\vee 0$, which  is precisely the tree obtained by the adjunction of a leaf starting from the (middle
of) the root. Then, the second component  $t\m 1$ is  
$ t^{(0)}\vee \cdots \vee (t^{(k)}+0)\vee 0 =  t^{(0)}\vee \cdots \vee t^{(k)}\vee 0 $, which  is precisely the tree
obtained by the adjunction of a leaf starting from the lowest vertex. Finally the first component  
$t\l 1$ is $ t^{(0)}\vee \cdots \vee (t^{(k)}+1)$ and here we use induction (on the degree of the last piece, that is
$t^{(k)}\ $) to prove that we get the union of the trees obtained from $t$ by adjoining a leaf to the other vertices
and edges.

$\arbreneufsixun $

$\arbreneufsixdeux $
\M

Let us start with a tree $s\in T_{n+1}$.  By deleting the most right  leaf we get a tree $t\in T_{n}$ such that $s$ 
 belongs to $t+1$ by the preceding argument. This proves existence. On the other hand, if $t\in T_n$ is such that
$s$ belongs to $t+1$, then obviously by deleting the most right  leaf  of $s$ we recover $t$. This proves unicity. \hfill
$\square$
\M

\M

\N {\bf 9.7. Remark.} The proof of $ {\underline n}+{\underline 1}={\underline {n+1}}$ in the binary case as given in
Proposition 2.3 was using the definition of the addition in terms of the poset structure of $Y_n$. Here the proof of the
analogue result is done by using the recursive definition of the sum. Since, in the binary case, one can also define the
addition through a recursive formula, this proof works also under the following
modification: the trees in
$s+1$ are obtained from $s$ by adding a new leaf starting from the middle of the edges on the right side. So all these
trees are still binary.
\M

\N {\bf 9.8. Involution.} Observe that for a planar tree {\it symmetry} around the axis passing
through the root defines an involution $\ss$ on $T_n$ and therefore also on
$\TT_n$. For instance $\ss(\ \arbreA) = \arbreA$, $\ss (\ \arbreB) = \arbreC$ and $\ss(\ \arbredeuxdeux)=
\arbredeuxdeux$.

By induction it comes immediately
 $$\eqalign{
\ss(x\l y) &= \ss (y)\r\ss (x)\cr
\ss(x\r y) &= \ss (y)\l\ss (x)\cr
\ss(x\m y) &= \ss (y)\m\ss (x)\cr
}$$
and therefore one gets:
$$\ss(x+y) = \ss(y) +\ss(x).$$
\M

\N{\bf 9.9. Filtration.} Define a filtration on $\TT_n$ by
$$ F^i \TT_n := \TT_{n,1}\ \cup \cdots \ \cup \ \TT_{n,i}.$$
So we have
$$\emptyset = F^0\TT_n \subset  \YY_n = F^1\TT_n\subset \cdots \subset  F^i\TT_n\subset \cdots  \subset  F^n\TT_n =
\TT_n.$$
\N {\bf Claim.} {\it If $x\in  F^i\TT_p$ and $y\in  F^j\TT_q$, then $x+y\in  F^{{\rm sup}(i,j)}\TT_{p+q}$.}
\M

The proof is straightforward by induction. In fact if $x\in  \TT_{p,i}$ and $y\in  \TT_{q,j}$, then $x+y\in 
\TT_{p+q,{\rm sup}(i,j) }\ \cup \ \cdots \ \cup \TT_{p+q,i+j}\ $.

In particular, fix $i$ and consider the sets $\TT_{n,i}$ for $n\geq i$. Then there is a well-defined sum
$$ \TT_{p,i}\times  \TT_{q,i}\longrightarrow  \TT_{p+q,i}\ , \ (x,y) \mapsto x\  {+_i}\  y$$
which consists in taking in $x+y$ only the trees belonging to  $ \TT_{p+q,i}$. From the preceding claim the operation 
${+_i}$ is associative. Indeed, in $x+y+z$ the elements in $\TT_{p+q+r,i}$ cannot come from elements in 
$ \TT_{p,j}\times  \TT_{q,i} \times  \TT_{r,i}$ with $j>i$, for instance, since ${\rm sup}(j,i,i) = j>i$.

For $i=1$ one has $ \TT_{n,1} = \YY_n$ and one recovers the sum of planar binary trees devised in 2.1.
\B

Summarizing what we have proved until now, we get the following 
\M

\N {\bf 9.10. Corollary.} {\it The set  $\TT_{\infty} := \cup_{n\geq 0}\, \TT_n$ of
groves is an involutive  graded monoid for $+$. The maps 
$$\displaylines{
\NN\longrightarrow \TT \longrightarrow \NN\cr
n\mapsto \u n \quad , \quad t\mapsto {\rm deg}\  t\cr
}$$
are morphisms of monoids.

The sets  $\TT_{\infty , i} := \TT_0 \cup \bigcup_{n\geq i}  \TT_{n,i}$ is an involutive  graded monoid for $+_i$.}
\BB

\N {\bf 10. Permutation-like notation of trees.} We extend the permutation-like notation introduced in section 3 to all
planar trees. 
\M

\N {\bf 10.1. Definition.} By definition the {\it name} of the unique element of $T_0$ is
$0$, and the name $w(x)$ of the element $x\in T_n$,
$n\geq 1$, is a finite sequence of strictly positive integers obtained inductively as follows. 

If $x = x^{(0)}\vee x^{(1)}\vee \cdots \vee x^{(k)}\in
 T_n$, then 
$$w(x):= (w(x^{(0)}), n,w(x^{(1)}), n,\cdots , n , w(x^{(k)}))$$
with the convention that we do not write the zeros. If there is no possibility of confusion we simply write 
$w(x):= w(x^{(0)}) \ n\cdots  n\   w(x^{(k)})$ (concatenation). Observe that, except for $0$, such a sequence is made of
$n$ integers and the integer $n$ appears $k$ times.  The name of the unique element of $T_1$ is therefore $1$.
\M

\N {\bf 10.2. Bijection with the planar trees.} When an element in $T_n$ corresponds to a planar tree $x$ and to a
name $w(x)$, we will say that $w(x)$ is the name of the tree $x$. In low dimension we get for instance

$$\matrix{ x &= &\vert   & \arbreA& \arbreB & \arbreC &\arbredeuxdeux 
&\arbreun & \arbreuntroistrois & \arbretroistroistrois\cr
w(x) &= &0& 1&12 & 21 & 22 & 123 & 133 & 333 \cr
}$$

Recall that symmetry around the root axis induces an involution on $T_n$. If $a_1 \cdots a_n$ is the name of the tree
$x$, then the name of $\ss(x)$ is $a_n \cdots a_1$.
\M

\N{\bf 10.3. Test for sequences.} Given a sequence of integers, is it the name of a
planar treeÊ? First, check that the largest integers in the sequence are the length of
the sequence. Second,  check that the maximal subsequences not containing the
largest integers are names of trees. Example: 14218812 is the name of a tree,
323 is not.
\M

\N {\bf  10.4. Relationship with the symmetric group.} In 3.5 we mentioned the relationship between the permutations
and the planar binary trees. There is a similar relationship between ordered partitions of $\{ 1, 2, \cdots , n\}$ and
$T_n$.
\M

\N {\bf 10.5. Tricks for computation.} For some special trees the computation of the
Left or Right or Middle sum is easy. First recall that  
$$-\vee - \cdots -\vee -  = - {\tt max} - \cdots - {\tt max} -\ ,$$
 where {\tt max} stands
 for the largest integer (i.e the length) of the word to which it belongs. Then it is easy to check that
$$\eqalign{
(- {\tt max}- \cdots - {\tt max}) \l - &=- {\tt max}- \cdots - {\tt max} - \cr 
-  \r ({\tt max} - \cdots - {\tt max}-) &=-  {\tt max} - \cdots - {\tt max}-  \cr 
(- {\tt max}- \cdot - {\tt max}) \m (- {\tt max}- \cdot - {\tt max} - ) &= - {\tt max}- \cdot - {\tt max} - {\tt
max}- \cdot - {\tt max} -)  \cr  
(- {\tt max}- \cdot - {\tt max}-) \m ({\tt max}- \cdot - {\tt max} - ) &= - {\tt max}- \cdot - {\tt max} - {\tt
max}- \cdot - {\tt max} -)  \cr  
}$$
For instance:
$313 \l 12 = 51512, \quad 313 \m 12 = 51515.$
\B

\N {\bf 11. Polynomial algebra with tree exponents.} In this section we recall some results on dendriform trialgebras
from [LR3, LR4] and we deduce algebraic relations between the operations on $\TT_{\infty}$.
\M

\N {\bf 11.1. Dendriform trialgebras.} By definition a {\it dendriform trialgebra} is a $K$-vector
space  $A$ equipped with three binary operations
$$
\eqalign{
\prec &: A\otimes A\to A,\cr
\succ &: A\otimes A\to A,\cr
\cdot &: A\otimes A\to A,\cr
}$$
which satisfy the following axioms:
$$\left\{\eqalign {
(a \prec  b) \prec  c &= a \prec  (b * c)\ , \cr
(a \succ  b) \prec  c &= a \succ  (b \prec  c)\ , \cr
(a * b) \succ  c &= a \succ  (b \succ  c)\ , \cr
\cr
(a \succ  b) \cdot  c &= a \succ  (b \cdot  c)\ , \cr
(a \prec  b) \cdot  c &= a \cdot  (b \succ  c)\ , \cr
(a \cdot  b) \prec  c &= a \cdot  (b \prec  c)\ , \cr
\cr
(a \cdot  b) \cdot  c &= a \cdot  (b \cdot  c)\ . \cr
}\right.$$
for any elements $a, b$ and $c$ in $A$, with the notation
$$a*b := a\prec b + a\succ b + a\cdot b\  .$$
Adding up all  the relations shows that the operation $* $ is associative.
\M

\N {\bf 11.2. A dendriform trialgebra associated to planar trees.} Let  $K[T_{\infty}']$ be the vector
space generated by the elements
$X^x, $ for $x\in T_n\  (n\geq 1)$ (we exclude $1\in T_0$ for a while). For a union of trees $\oo = \cup_i\, x_i$, we
introduce the notation
$$ X^{\oo} := \sum_i X^{x_i}\ . $$
We define three operations on  $K[T_{\infty}']$ by the formulas
$$\eqalign{
X^x\prec X^y &:= X^{x\l y}\, , \cr
X^x\succ X^y &:= X^{x\r y}\, ,  \cr
X^x\cdot X^y &:=  X^{x\m y}\,  ,\cr
}$$
and distributivity.

We now recall the following main result:
\M

\N {\bf 11.3. Theorem (Universal property)} [LR3, LR4]. {\it The vector space $K[T_{\infty}']$ 
equipped with the three operations $\prec\  ,\  \succ$ and $\cdot$ as above is a dendriform trialgebra.
Moreover it is the free dendriform trialgebra on one generator (namely $X$). \hfill $\square$}
\M
Therefore $K[T_{\infty}']$ equipped with $*$ is an associative algebra. We make it into an associative and unital algebra 
 $K[T_{\infty}]= K\oplus K[T_{\infty}']$ by adding the vector space $K$, generated by $1= X^0$.
\M

\N {\bf 11.4. Relations in $\TT_{\infty}$.} The dictionnary comparing the operations on trees and on the free dendriform
algebra is the following:
\B
 {\vbox{
\hrule
\halign{&\vrule#&\strut # \cr
\ $\TT_{\infty}$	&& &&\  $\cup$  &&\  $\l$ &&\  $\r$ &&\  $\m$ &&\  + && \cr
\noalign{\hrule}
\ $K[T_{\infty}]$	&& &&\  $+$  &&\  $\prec$ &&\  $\succ$ &&\  $\cdot$ &&\  $*$ && \cr
\noalign{\hrule}
  }} 
\B
By theorem 11.3 associativity of the operation $*$ implies associativity of the operation $+$ on planar trees (this proves
Proposition 9.4).

Moreover, the relations in the dendriform algebra give for $x,y,z \in \TT_{\infty}$:

$$\left\{\eqalign {
(x \l  y) \l  z &= x \l  (y + z)\ , \cr
(x \r  y) \l  z &= x \r  (y \l  z)\ , \cr
(x + y) \r  z &= x \r  (y \r  z)\ , \cr
\cr
(x \r  y) \m  z &= x \r  (y \m  z)\ , \cr
(x \l  y) \m  z &= x \m  (y \r  z)\ , \cr
(x \m  y) \l  z &= x \m  (y \l  z)\ , \cr
\cr
(x \m  y) \m  z &= x \m  (y \m  z)\ . \cr
}\right. \eqno (*)$$
\M

\N {\bf 11.5. Corollary.} {\it Any planar tree $x\in T_n$ can be written as a composite (using $\l\ , \r \ , \m $) of
$n$ copies of $1\in T_1$ with a suitable parenthesizing. This universal expression of $x$, denoted $w_x(1)$, is unique
modulo the relations $(*)$. \hfill $\square$}
\M
\N Examples: $w_{21}(1) = 1\l1,\  w_{12}(1) = 1\r 1,\  w_{22}(1) = 1\m 1, $

 \N $w_{122}(1) = 1\r 1\m 1,\  w_{322}(1) = 1\l (1\m1) .$
\M

\N {\bf 11.6. Remarks about notation.} In [LR3, LR4] the linear generators are denoted $x$
instead of $X^x$. We adopt this different  notation here to avoid confusion with the operations in
$\TT_{\infty}$.

In [LR3, LR4] the symbols $\l$, $\r$ and $\m$ are used to denote operations of an associative
trialgebra which is the Koszul dual structure of dendriform trialgebra. We have given them a
completely different meaning here.
\B

\N {\bf 12. Multiplication.} Since the associative algebra $K[T_{\infty}']$ is free on one generator when
considered as a dendriform trialgebra, one can perform composition of polynomials with tree
exponents. Though the composite of monomials, that is $(X^x)^y$ where $x$ and $y$ are
planar trees, is not a monomial, it turns out that it is $X$ to the power of some grove. Hence one
can define the multiplication of planar trees as being this grove and then extend this
multiplication to any groves.
\M
\N {\bf 12.1. Definition.} Let   $x$ and $y$ be planar trees. By definition the {\it product}
$x\times y$ is
$$ x \times y := w_x(y).$$
where $ w_x(1)$ is the universal expression of $x$ (cf. Corollary 11.5). In other words we replace
all the copies of $1$ by copies of $y$ in this universal expression.

Observe that the above definition of the product has a meaning even when $y$ is a grove since the 
Right sum, the Left sum and the Middle sum of groves are well-defined.  We extend the multiplication to $x$ being a
grove by distributivity on the left with respect to disjoint union:
$$(x \cup x')\times y = w_x(y)\cup w_{x'}(y).$$
So we have defined the product of two groves. This product is a grove since it is obtained by the
operations $\l$, $\r$ and $\m$. Observe that the product is not commutative.
\M

 \N {\bf 12.2. Proposition.} {\it The multiplication $\times$ on groves is
distributive  on the left with respect to the Left sum, the Right sum, the Middle sum  and the sum (but not on the right).}
\M
\N {\it Proof.} The formula $w_{x+x'}(1) = w_x(1)+ w_{x'}(1)$ follows from 
$$\eqalign{
w_{x\l x'}(1) &=  w_x(1)\l w_{x'}(1)\ ,\cr
w_{x\r x'}(1) &=  w_x(1)\r w_{x'}(1)\ ,\cr
w_{x\m x'}(1) &=  w_x(1)\m w_{x'}(1)\ .\cr
}$$
These three formulas follow inductively from the properties of the function $w$.\hfill $\square$
\M

 \N {\bf 12.3. Proposition.} {\it
The multiplication of groves is associative with neutral element on both sides the tree
$\arbreA = 1$.} 
\M
\N {\it Proof.}  Interpreted in terms of dendriform trialgebra, the multiplication of planar trees is
composition of monomials. Since composition is associative, the multiplication of planar trees is
associative.

Since $w_1(1)=1$, we get $w_1(y)=y$ and so $1\times y = y$. On the other side $x\times1 = x$
is a tautology. \hfill $\square$
\M

\N {\bf 12.4. Theorem.} {\it With the notation ${\u  n} = \cup_{x\in T_n}\ 
x$, one has 
$${\u  n} \times {\u  m} = {\u  {nm}}.$$}
\N {\it Proof.} Since, by Proposition 12.2 the multiplication is distributive on the left with respect to addition, we
get 
$$\eqalign{
{\u  n} \times {\u  m} &= ({\u  1}+ \cdots + {\u  1})\times {\u  m} = {\u  1}\times {\u  m}
+ \cdots + {\u  1}\times {\u  m}\cr
&= {\u  m} + \cdots + {\u  m}=  {\u  {nm}}.\cr
}$$
 { }\hfill $\square$
\S

\N {\bf 12.5. Proposition (recursive property).} {\it Let $x= x^{(0)}\vee \cdots \vee x^{(k)}$ be a planar tree and let $y$
be a grove. The multiplication is given recursively by the formulas
$$x\times y =
 (x^{(0)}\times y ) \r y \m  (x^{(1)}\times y )\m y \cdots \m y \l  (x^{(k)}\times y)$$
and $0\times y = 0$. }
\M
\N {\it Proof.} First, observe that because of the relations (*) there is no need for parenthesis in this formula. It suffices
to show the equality
$$w_x(y) =
 w_{x^{(0)}}(y ) \r y \m  w_{x^{(1)}}( y )\m y \cdots \m y \l w_{x^{(k)}}(y)\ ,$$
which is a consequence of
$$x =  x^{(0)}\vee \cdots \vee x^{(k)}=
  x^{(0)} \r 1 \m   x^{(1)} \m 1 \cdots \m 1\l  x^{(k)} \ .
$$
This last formula is easily proved by using the definition of the three operations.
\hfill $\square$
\M
\N {\bf 12.6. Proposition (involution).} {\it For any groves $x$ and $y$ one has
$$\ss(x\times y) = \ss (x) \times \ss (y).$$}
\N {\it Proof.} It is a consequence of the previous Proposition and of the relations between the involution and the
three operations (cf. 9.8).
\hfill $\square$
\M

\N {\bf 12.7. Summary.} On the set of groves $\TT_{\infty} =
\bigcup_{n\ge 0}\TT_n$ there are defined operations + and $\times$ such that

$\bullet$ the addition + is associative, distributive both sides with respect
to $\cup$, with neutral element $0=\vert$ , but is not commutative,

$\bullet$ the multiplication $\times$ is associative, distributive on the left
with respect to the sum + and to the disjoint union $\cup$ (but not right distributive), with neutral
element (both sides) $1=\arbreA$, but is not commutative,

$\bullet$ the involution $\ss$ on $\TT_{\infty}$ satisfies $\ss(x+y) = \ss (y) + \ss (x)$ and 
$\ss(x\times y) = \ss (x) \times \ss (y)$,

$\bullet$  the maps $\NN \to \TT_{\infty}, n\mapsto {\underline n}=\bigcup_{x\in Y_n} x$ and 
 $\deg : \TT_{\infty} \to \NN$
(degree) are compatible with $+$ and $\times$. The composite is the identity of
$\NN$.

$\bullet$ The quotient map $\TT_{\infty} \to \YY_{\infty}$, which consists in forgetting about the planar trees which
are not binary, is both additive and multiplicative.
\M

\N {\bf 12.8. Question.} It would be interesting to know if one can also put an internal  multiplication on $\TT_{\infty , i}$
for
$i>1$, cf. 9.9.
\vfill
\eject

\N {\bf Appendix to part II: tables for planar trees.}
\B

\N {\bf II.A.1. Addition table.} Recall that $0$ is the neutral element for $+$ , so
$$0+x = x = x+ 0.$$
In the following table we omit the $\cup$ sign. The first line is $x\r y$,  the
second line  is $x\l y$ and the
third line  is $x\m y$.
\B

{\petit{
{\vbox{
\hrule
\halign{&\vrule#&\strut \hfil# \cr
$\ x + y$	&&\ 1 		&&\ 12 			&&\ 21		&&\  22  &&\cr
\noalign{\hrule}
\ 1  	&&\ 12 		&&\ 123\  213\ 223			&&\ 131		&&\ 133 && \cr 
	&&\ 21     	&&\ 312		&&\ 321		&&\ 322 && \cr
	&&\ 22     	&&\ 313		&&\ 331		&&\ 333 && \cr
\noalign{\hrule}
\ 12 	&&\ 123	    	&&\  1234\ 1314\ 1334 		&&\ 1241		&&\ 1244 && \cr
	&&\ 131	&&\ 1412			&&\ 1421 	 	&&\ 1422 &&\cr
	&&\ 133	&&\ 1414			&&\ 1441 	 	&&\ 1444 &&\cr
\noalign{\hrule}
\ 21 	&&\ 213 &&\  2134\ 3124\ 3214\ 3224\ 3134 &&\ 2141 	&&\ 2144 && \cr 
	&&\  312	\ 321\ 322		&&\ 4123\  4213\  4312\ 4223\ 4313 		&&\  4131\ 4321\ 4331		&&\ 4133\ 4322\ 4333 &&
\cr
	&&\  313	&&\ 4124\  4214\ 4224		&&\  4141	&&\ 4144 && \cr
\noalign{\hrule}
\ 22	&&\ 223		&&\  2234\ 3314\ 3334		&&\ 2241		&&\ 2244 && \cr
	&&\ 331 		&&\ 4412		&&\ 4421		&&\ 4422 && \cr
	&&\ 333 		&&\ 4414		&&\ 4441		&&\ 4444  && \cr
\noalign{\hrule}
  }}}
}
\BB
\normal 

\N {\bf II.A.2. Mutiplication table.} Recall that 1 is the neutral element for $\times$ , so 
$$1 \times x =x =x\times 1.$$
\B
\petit{
{\vbox{
\hrule
\halign{&\vrule#&\strut \hfil# \cr
$\ x \times y$	&&\ 12 &&\ 21 &&\ 22 &&\ 12 \ 21	 && \  123 &&  \  133 &&\cr
\noalign{\hrule}
\ 12  &&\ 1234	&&\ 2141 &&\ 2244 &&\ 1234\ 1314\ 1241\ 1334   &&\ 125126\ 124156 &&\ 134166 &&\cr 
\  	  &&\ 1314   &&\ 		   &&\             &&\ 2134\ 3124\ 3214\ 3224	&&\ 123456\ 124456  &&\ 133466 &&\cr
\  	  &&\ 1334	 &&\ 		   &&\             &&\ 2141\ 3134	                   &&\ 125156 &&\ 134466 &&\cr
\noalign{\hrule}
\ 21  &&\ 1412	&&\ 4131 &&\ 4422 &&\ 1412\ 4123\ 4213\ 4223 &&\ 126123 &&\ 166133 &&\cr
\  	  &&\ 		       &&\ 4321 &&\          &&\ 4312\ 1421\ 4131\ 4313 &&    &&\  &&\cr
\  	  &&\ 		       &&\ 4331 &&\           &&\ 4321\ 4331		                  &&    &&\ &&\cr
\noalign{\hrule}
\ 22  &&\ 1414	&&\ 4141 &&\ 4444 &&\ 1414\ 1441\ 4124\ 4214 &&\ 126126 && \ 166166 &&\cr
\  	  &&\ 		      &&           	&&\          &&\ 4224\ 4141                     &&       &&\  &&\cr
\noalign{\hrule}
  }}
}
\vfill\eject

\normal
\N {\bf III. Final comments}
\M
$\bullet$ It is sometimes helpful to index chain complexes, not by the integers or even pair of integers (like in a
bicomplex), but by trees. Examples and the simplicial properties of the planar binary trees have been investigated by A.
Frabetti in [Fr]. Similarly, many small categories in algebraic topology have the natural numbers as objects (for instance
$\Delta, \Gamma$). In the work of Andr\'e Joyal [J] appears a category $\Theta$ whose objects are the planar trees as
considered here.
\M

$\bullet$ The relationship of dendriform algebras and trialgebras with other types of algebras like associative
dialgebras, associative trialgebras, Leibniz algebras, Zinbiel algebras, associative algebras has been treated in
[L2] and [LR4] in terms of operads. See also [Ch].
\M

$\bullet$ The operad of associative algebras comes from a set-operad (sometimes called the Barratt-Eccles operad) by
the functor which associates to a set the vector space based on it. The operad of dendriform dialgebras (resp. dendriform
trialgebras) does not come from a set operad since the sum of two trees is not a tree. However since the sum of two
groves is a grove, it is very close to being a set operad.
\M

$\bullet$ The associative algebra $ {K[Y_{\infty}]}$ has a dendriform structure, 
but has also a Hopf algebra structure (cf. [LR1]). Moreover these two structures are compatible, as was discovered by M.
Ronco [R1]. It turns out that $ {K[Y_{\infty}]}$ is the universal
enveloping dendriform dialgebra of the free {\it brace algebra} on one generator as proved by M. Ronco
in [R2, R3]. A similar result holds for the Hopf algebra $K[T_{\infty}]$.
\M

$\bullet$ Planar binary trees have been used by
Christian Brouder [Br] in place of natural numbers in order to index series which are
solutions of some differential equations of quantum field theory (the Schwinger-Dyson equations). The {\it
renormalization} of quantum electrodynamics is governed by a certain non-commutative and non-cocommutative Hopf
algebra, cf. [BF]. It turns out that this Hopf algebra is isomorphic to $K[Y_{\infty}]$, cf. [Fo], [H].
\BB

\centerline {\bf References}
\B

\N [Br] Ch. Brouder, {\it On the trees of quantum fields}, Eur. Phys. J. {\bf C 12} ({\oldstyle 2000}), 535--549.
\S
\N [BF] Ch. Brouder and A. Frabetti, ``Renormalization of QED with trees", Eur. Phys. J.   {\bf C 19}
({\oldstyle 2001}), 715--741. 
\S
\N [Ch] F. Chapoton, {\it Alg\`ebres de Hopf des permutah\`edres, associah\`edres et hypercubes. } Adv. Math. 150
({\oldstyle 2000}), no. 2, 264--275. 
\S
\N [Co] J.H.C. Conway, ``On Numbers and Games". Second edition. A K Peters, Ltd., Natick, MA, {\oldstyle 2001}.
\S
\N [Fo] L. Foissy, {\it Les alg\`ebres de Hopf des arbres enracin\'es d\'ecor\'es.} Th\`ese, Reims,  {\oldstyle 2001}.
\S
\N [Fr] A. Frabetti, {\it Simplicial properties of the set of planar binary trees}. J. Algebraic Combin. 13
({\oldstyle 2001}), no. 1, 41--65.  
\S
\N [H] R. Holtkamp, {\it Comparison of Hopf algebra structures on trees.} Preprint Bochum,  {\oldstyle 2001}.
\S
\N [J] A. Joyal, {\it Disks, duality and $\Theta$-categories}, preprint 1997, 6 pages.
\S
\N [L1] J.-L. Loday,  {\it Alg\`ebres ayant deux op\'erations associatives (dig\`ebres).} C. R. Acad. Sci. Paris S\'er. I Math. 
{\bf 321} ({\oldstyle 1995}), no. 2, 141--146.
\S
\N [L2] J.-L. Loday,  {\it Dialgebras},  in ``Dialgebras and related operads", Springer Lecture
Notes in Math. {\bf 1763} ({\oldstyle 2001}), 7--66.
\S
\N [LR1] J.-L. Loday,  and M. O. Ronco,   {\it Hopf algebra of the planar binary trees},  Adv. Math.
{\bf 139} ({\oldstyle 1998}), no. 2, 293--309.
\S
\N  [LR2]  J.-L. Loday, and M.O. Ronco, {\it  Order structure on the algebra of
permutations and of planar binary trees}. J. Alg. Comb. ({\oldstyle 2001}), to appear.
\S
\N  [LR3]  J.-L. Loday, and M.O. Ronco, {\it  Une dualit\'e entre simplexes standards et polytopes de Stasheff}. C. R. Acad.
Sci. Paris {\bf 333} ({\oldstyle 2001}), 81--86.
\S
\N  [LR4]  J.-L. Loday, and M.O. Ronco, {\it  Trialgebras and families of polytopes}. Preprint ({\oldstyle 2001}).
\S
\N  [R1]  M.O. Ronco, {\it  Primitive elements in a free dendriform algebra.} New trends in Hopf algebra theory (La
Falda, 1999), 245--263, Contemp. Math., {\bf 267}, Amer. Math. Soc., Providence, RI, {\oldstyle 2000}. 
\S
\N  [R2]  M.O. Ronco, {\it A Milnor-Moore theorem for dendriform Hopf algebras.} C. R. Acad. Sci. Paris S\'er. I Math.
{\bf 332} ({\oldstyle 2001}), no. 2, 109--114.
\S
\N  [R3]  M.O. Ronco, {\it  A Milnor-Moore theorem for some non-cocommutative
Hopf algebras}. Preprint ({\oldstyle 2000}).
\S
\N [St] R.P. Stanley,  ``Enumerative combinatorics". Vol. I. The Wadsworth and
Brooks/Cole Mathematics Series, {\oldstyle 1986}.

\B
 Institut de Recherche Math\'ematique Avanc\'ee

    CNRS et Universit\'e Louis Pasteur

    7 rue R. Descartes

    67084 Strasbourg Cedex, France (EU)

    Courriel : loday@math.u-strasbg.fr
\B
{\tt [arithmetree] } \hfill 27 novembre 2001

\end